\documentclass[11pt,a4paper,leqno]{article}
\usepackage{amsmath,amsthm,amsfonts}
\newtheorem{defin}{Definition}[section]
\newtheorem{lemma}[defin]{Lemma}
\newtheorem{propos}[defin]{Proposition}
\newtheorem{theor}[defin]{Theorem}
\newtheorem{corol}[defin]{Corollary}
\newenvironment{definition}{\begin{defin} \rm}{\hspace{2mm}\BOX \end{defin}}
\newenvironment{definition*}{\begin{defin} \rm}{\end{defin}}
\newcommand{\smallsep}{\vspace{-2.9mm}\itemsep -1.6mm}

\newcommand{\bd}{\begin{definition}}
\newcommand{\ed}{\end{definition}}
\newcommand{\bdefn}{\begin{definition*}}
\newcommand{\edefn}{\end{definition*}}
\newcommand{\bl}{\begin{lemma}}
\newcommand{\el}{\end{lemma}}
\newcommand{\bp}{\begin{propos}}
\newcommand{\ep}{\end{propos}}
\newcommand{\bt}{\begin{theor}}
\newcommand{\et}{\end{theor}}
\newcommand{\bc}{\begin{corol}}
\newcommand{\ec}{\end{corol}}
\newcommand{\bds}{\begin{displaystyle}}
\newcommand{\eds}{\end{displaystyle}}
\newcommand{\BOX}{\rule{2mm}{2mm}}
\newcommand{\Z}{${\mathbb Z}_{2}$}

\newcommand{\CCq}{\mathbb C(q)}
\newcommand{\FF}{{\mathbb F\hskip.5pt}}
\newcommand{\KK}{{\mathbb K\hskip.5pt}}
\newcommand{\groupring}
{\mbox{${\mathbb C} \hspace{-0.3mm}\cdot\hspace{-0.4mm} S_{n}$\ }}
\newcommand{\Sa}{G_{n}}
\newcommand{\qSa}{G_{n}(q)}
\newcommand{\qSprime}{G_{n}^{\,\prime}(q)}
\newcommand{\qSdash}{G_{n}^{\hspace{0.65mm}\text{\bf-}}(q)}
\newcommand{\qSast}{G_{n}^{\hspace{0.15mm}\ast}(q)}
\newcommand{\ASA}{Se_{n}(q)}
\newcommand{\loc}{Se_{n}^{\hspace{0.15mm}\circ}(q)}
\newcommand{\CS}{{\mathbb C}(q)[X]}
\newcommand{\xseq}{(x_{1},\ldots,x_{n})}
\newcommand{\qseq}{(q_{\qs 1},\ldots,q_{\qs n})}
\newcommand{\qs}{\hspace{0.3mm}}
\newcommand{\stp}{\begin{picture}(59,20)
\put(10,5){\mbox{\begin{scriptsize}$>$\hspace{-0.8mm}$-$\end{scriptsize}}}
\end{picture}}
\newcommand{\stpscript}{\begin{picture}(45,18)
\put(7,4){\mbox{\begin{tiny}$>$\hspace{-0.7mm}$-$\end{tiny}}}
\end{picture}}
\begin{document}
\thispagestyle{empty}
\begin{center}
%
%
{\LARGE 
Affine Sergeev Algebra and\\
$q\hskip1pt$-Analogues of the Young Symmetrizers for\\[2mm]
Projective Representations of the Symmetric Group}
\\[6mm]
%
%
{\Large A.\,R.\,Jones, M.\,L.\,Nazarov}\\[6mm]
%
%
\parbox{5.42in}
{We study a $q$-deformation for the semi-direct product of the symmetric
group $S_{n}$ with the Clifford algebra on $n$ anticommuting generators.
We obtain a $q$-version of the projective analogue
for the classical Young symmetrizer introduced by the second author.
Our main tool is an analogue of the Hecke algebra of complex valued functions
on the group $GL_{n}$ over a $p$-adic field relative to the Iwahori subgroup.}
\end{center}
\medskip
\section{Introduction}
\setcounter{equation}{0}
The classical representation theory of the symmetric group $S_{n}$ is
well known.
The irreducible representations of $S_{n}$ over the complex field
${\mathbb C}$ are labelled by partitions $\omega$ of the integer~$n$
and realised through several constructions \cite{Jam}.
One of the most useful realisations employs primitive idempotents in the
group ring \groupring now known as Young symmetrizers \cite{Young}.
These idempotents in \groupring correspond to the
Young standard tableaux $\Omega$ with $n$ entries \cite{Mac}.
Suppose that a tableau $\Omega$ has shape $\omega$.
Let $S_\Omega$ and $S^{\hskip1pt\ast}_\Omega$ be the subgroups in $S_n$
\mbox{respectively} preserving the rows and columns of $\Omega$ as sets.
The {\it Young symmetrizer} corresponding to $\Omega$ is
\begin{equation}
\label{young_symmetrizer}
\frac{n_\omega}{n!}\,
\sum_{s\in S_\Omega}\,\sum_{t\in S^{\hskip1pt\ast}_\Omega}\,
s\,t\cdot\operatorname{sgn}(t)
\end{equation}
where $n_{\omega}$ denotes the number of all standard tableaux with shape
$\omega$. We will write $\Omega(i,j)$ for the entry of the tableau $\Omega$
in the row $i$ and column $j$.

The $q$-analogues of Young symmetrizers are also known\qs;
for instance, see \cite{Gyo}.
Here the group ring \groupring is replaced by the
{\it Hecke algebra} $H_{n}(q)$ of type $A_{n-1}$.
This is the algebra over the field $\CCq$ generated by the elements
$T_1,\ldots,T_{n-1}$ with the relations
\begin{equation}
\label{sergeev_def1}
\begin{array}{rcl}
( T_{k} - q )( T_{k} + q^{-1} ) &=& 0\,; \\
T_{k} T_{k+1} T_{k} &=& T_{k+1} T_{k} T_{k+1}\,; \\
T_{k} T_{l} &=& T_{l} T_{k}\,, \hspace{4mm} l\neq k,k+1\,.
\end{array}
\end{equation}
As usual, we will write $T_s=T_{k_1}\cdots \,T_{k_l}$ for any choice of
reduced decomposition $s=s_{k_1}\cdots \,s_{k_l}$ in the
standard generators $s_{1},\ldots,s_{n-1}$ of the symmetric group $S_n$.
For each partition $\omega$ there are two distinguished standard tableaux  
of shape $\omega$ obtained by inserting the symbols $1,2,\ldots,n$ into
the Young diagram of shape $\omega$ by rows and columns.
These are called the {\it row} and {\it column tableaux\/} and
are denoted by $\Omega^{r}$ and $\Omega^c$ respectively.
Consider the elements of the algebra $H_n(q)$
$$
E=\sum_{s\in S_{\Omega^r}}\,T_s\cdot q^{\hskip1pt l}
\,,\ \quad
E^{\hskip1pt\ast}=\sum_{s\in S^{\hskip1pt\ast}_{\Omega^c}}\,T_s\cdot(-q)^{-l}
$$
where $l$ is the length of the reduced decomposition of $s$.
Now choose the element $s\in S_n$ such that
$s(\hskip1pt\Omega^r(i,j))=\Omega^c(i,j)$.
Up to a certain non-zero multiplier from $\CCq$, the product
$E_\omega=E\,T_s^{\hskip1pt-1}E^{\hskip1pt\ast}T_s\in H_n(q)$
is the $q$-analogue \cite[Section 2]{Gyo} of the Young symmetrizer 
(\ref{young_symmetrizer}) corresponding to $\Omega=\Omega^r$.
That $q$-analogue is a primitive idempotent in the algebra $H_n(q)$. 

Another approach to these $q$-analogues has been proposed by
Cherednik \cite{Cher2}\qs; it employs the {\it affine Hecke algebra} $He_n(q)$.
This is the algebra over $\CCq$ generated by $T_1,\ldots,T_{n-1}$ and the
pairwise commuting invertible elements $Y_1,\ldots,Y_n$ subject to the
relations
$$
T_{k}Y_{k}T_k\;=\;Y_{k+1}\,;
\qquad
T_{k}Y_{l}\;=\;Y_{l}T_{k}\,,\quad l \neq k,k+1
$$
along with the relations (\ref{sergeev_def1}).
Let $\chi$ be a character of the subalgebra in $He_n(q)$
formed by all Laurent polynomials in $Y_1,\ldots,Y_n$.
This subalgebra 
is maximal commutative\qs; see \cite{Cher4}.
Denote by $M_\chi$ the representation of $He_n(q)$
induced from the character $\chi$\,, it is called a  
{\it principal series representation} \cite{Rog}.
The vector space $M_{\chi}$ can be identified with the algebra
$H_{n}(q)$ so that the generators $T_{1},\ldots,T_{n-1}$ act via
left multiplication.
The action of the generators $Y_1,\ldots,Y_n$ in $M_\chi$ is
determined by $Y_k\cdot1=\chi(Y_k)$.
Note that the group $S_n$ acts on the characters $\chi$ in a
natural way: $w\cdot\chi(Y_k)=\chi(Y_{w^{-1}(k)})$ for any $w\in S_n$.
The element $E_{\omega}\in H_n(q)$
has appeared in \cite[Section 3]{Cher2} in the following guise: choose 
$w$ such that $s\,w=w_0$ is the element of the maximal length in $S_n$ 
and determine the character $\chi$ by
$$
w\cdot\chi(Y_k)\,=\,q^{\hskip1pt2(j-i)},
\ \quad
k = \Omega^{r}(i,j)\,.
$$
Then the operator of right multiplication in $H_{n}(q)$ by the element
$E_{\omega} T_w$ is an intertwining operator
$M_{w\cdot\chi} \rightarrow M_\chi$ over the algebra $He_{n}(q)$,
cf.\ \cite[Theorem 6.3]{Naz}.
This approach gives an explicit formula for the element $E_{\omega}$
different from that given above, cf.\ \cite[Theorem 5.6]{Naz}.

It was Schur \cite{Sch} who discovered non-trivial central
extensions of the \mbox{symmetric group $S_n$}.
In other words, the group $S_n$ admits
projective representations which cannot be reduced to linear ones. 
The analogues of the Young symmetrizers (\ref{young_symmetrizer}) for
these representations were constructed by the second author
in \cite{Naz} using the approach of \cite{Cher2}.
In \cite{Naz} the role of the group ring ${\mathbb C} \cdot S_n$
is taken by the crossed product $\Sa$ of the group $S_{n}$
with the Clifford algebra over $\mathbb{C}$ on $n$ anticommuting generators.
The irreducible $\Sa$-modules are parametrised by the partitions
$\lambda$ of $n$ with pairwise distinct parts.
Note that the algebra $\Sa$ has a natural \Z-grading and here the
notion of \Z-graded irreducibility is used.
These modules can be constructed from the projective representations
of the group $S_{n}$, see for instance \cite{Stem}.

A $q$-analogue of the algebra $\Sa$ was
introduced by Olshanski in \cite{Olsh} by generalizing the double commutant
theorem from \cite{Jim}.
It is the $\CCq$-algebra $\qSa$ generated by the
Hecke algebra $H_{n}(q)$ and the Clifford algebra with the generators
$C_1,\ldots,C_n$ subject to the relations 
(\ref{sergeev_def2}) and (\ref{sergeev_def3}) below.
The algebra $\qSa$ is \Z-graded so that
$\operatorname{deg}\,T_{k} = 0$ and $\operatorname{deg}\,C_{l} = 1$
for all possible $k$ and $l$.
The aim of the present paper is to construct the $q$-analogues in $\qSa$
of the projective Young symmetrizers from \cite{Naz}.
This is done by extending the approach of \cite{Cher2}.

In Section 2 we examine some basic properties of the algebra $\qSa$.
In particular, we show that the algebra $\qSa$ is semisimple\qs;
see Proposition \ref{semisimple}.
In Section 2 we also develop the combinatorics of shifted tableaux for
partitions of $n$ with pairwise distinct parts. 

In Section 3 we introduce the main underlying object of the present paper.
This is the algebra $\ASA$ over $\CCq$ generated by $\qSa$ and the
pairwise-commuting invertible elements $X_1,\ldots,X_n$
subject to the relations (\ref{affine_rels1}) and (\ref{affine_rels2}).
It is the analogue for $\qSa$ of the affine Hecke algebra $He_n(q)$.
We will call it the affine Sergeev algebra in honour of the author of
\cite{Serg} who extended the double commutant theorem of Schur and Weyl
from the group $S_n$ to the algebra $\Sa$.
The subalgebra in $\ASA$ formed by all Laurent polynomials in 
$X_1,\ldots,X_n$ is maximal commutative\qs; see Proposition \ref{max_comm}.
By studying intertwining operators between the principal series representations
of $\ASA$ relative to this subalgebra, we obtain for the algebra $\qSa$
an analogue of the element $E_\omega T_w\in H_{n}(q)$.

There exists a homomorphism $\iota : \ASA \rightarrow \qSa$ which is identical
on the subalgebra $\qSa \subset \ASA$\qs; see Proposition \ref{Murphy_hom}.
The images $J_{1},\ldots,J_{n}$ of the generators $X_{1},\ldots,X_{n}$ relative
to this homomorphism are called the Jucys-Murphy elements of the algebra
$\qSa$; cf.\ \cite{Cher2} and \cite{Jucy2}.
These elements play a major role in the present article.

In Section 4 we introduce a certain extension $\FF$ of the field $\CCq$.
This will be a splitting field for the
semisimple algebra \text{$\qSa$ over $\CCq$.}
As usual let us write for any non-negative integer $m$
$$
[m]_q=\frac{q^m-q^{-m}}{q-q^{-1}}\,.
$$
The field $\FF$ is obtained by adjoining to $\CCq$ a square root of
$[m]_{q^2}$ for each $m=2,\ldots, n$.
We assign to each standard shifted tableau
$\Lambda$ with $n$ entries an element $\psi_{\Lambda}$ in the extended algebra
$\qSa \otimes_{\,\CCq} \FF = \qSprime$.
This construction employs the fusion procedure introduced by Cherednik
\cite{Cher2}\qs; see \text{Theorem \ref{fusion_theorem}} of the present paper.
Each element $\psi_{\Lambda}$ is an eigenvector for the left multiplications
by the Jucys-Murphy elements $J_{1},\ldots,J_{n}$.
Moreover, for any $1 \leq k_{1} < \cdots < k_{p} \leq n$,
we have the equality in the algebra~$\qSprime$
$$
J_{k} \cdot C_{k_{1}}\cdots \,C_{k_{p}} \,\psi_{\Lambda}\,=\,
\Bigl(\,
[m+1]_{q^2}-[m]_{q^2}\mp(q-q^{-1})\sqrt{\,[m+1]_{q^2}[m]_{q^2}\!}\,\,
\Bigr)
\cdot C_{k_{1}}\cdots \,C_{k_{p}} \,\psi_{\Lambda}
$$
where $m=j-i$ for $k=\Lambda(i,j)$ 
and the sign $\mp$ depends on whether the index $k$
is contained in the set $\{ k_{1},\ldots,k_{p} \}$ or not.

Now let $\Lambda$ run through the set of all standard shifted tableaux
corresponding to a fixed partition $\lambda$ of $n$ with pairwise distinct
parts\qs;
here the number of parts will be denoted by $\ell_{\lambda}$.
Let us denote the row and column tableaux by $\Lambda^{r}$
and $\Lambda^{c}$, respectively.
Then the elements $C_{k_{1}} \cdots \,C_{k_{p}} \,\psi_{\Lambda}$ form a
basis in the left ideal $V_{\lambda} \subset \qSprime$ generated
by the element \mbox{$\psi_{\Lambda^{r}}$\qs;} see Theorem \ref{V_theorem}.
The vector space $V_{\lambda}$ over $\FF$ becomes a $\qSprime$-module
under left multiplication.
This module is not irreducible, but splits into a direct sum of
$2^{[\,\ell_{\lambda}/2\,]}$ copies of a certain irreducible
$\qSprime$-module $U_\lambda$\qs;
see Corollary \ref{splitting}.
Here we again mean the \Z-graded irreducibility.
The modules $U_\lambda$ are non-equivalent for distinct $\lambda$
and form a complete set of irreducible $\qSprime$-modules\qs;
see Corollary \ref{complete_set}.
Moreover, they remain irreducible on passing to any extension of the field
$\FF$\qs; see Theorem \ref{V_results}.

Define the element $w_{\Lambda^r}\in S_n$ by the equality
$s_{\Lambda^r}w_{\Lambda^r}=w_0$ where
$s_{\Lambda^r}(\Lambda^r(i,j))=\Lambda^c(i,j)$.
The product $\psi_{\Lambda^r}T_{w_{\Lambda^r}}^{\,-1}\!\in \qSprime$
is our analogue of the above described element $E_\omega\in H_n(q)$.
Corollaries \ref{row_div} and \ref{col_div3} justify this claim.

We are grateful to A.\,O.\,Morris and G.\,I.\,Olshanski
for their interest in our work.
We are also grateful to G.\,D.\,James and A.\,E.\,Zalesskii
for their kind advice. Our work was supported by the
Engineering and Physical Sciences Research Council.
It was also supported by the Isaac Newton Institute for
Mathematical Sciences at Cambridge. 

\section{The Algebra \boldmath{$\Sa$} and its \boldmath{$q$}-Deformation}
\setcounter{equation}{0}
Let us recall some results on the algebra $\Sa$
and examine its $q$-deformation $\qSa$ introduced by Olshanski \cite{Olsh}
where it was referred to as the {\it Hecke-Clifford superalgebra}.
The algebra $\Sa$ is defined as the crossed product
of the symmetric
group $S_{n}$ with the Clifford algebra over the field ${\mathbb C}$ on $n$
generators\qs;
it has was studied by Sergeev \cite{Serg}.
Here the Clifford algebra is generated over ${\mathbb C}$ by
elements $C_{1},\ldots,C_{n}$ obeying the relations
\begin{equation}
\label{sergeev_def2}
C_{k}^{2} \;=\; -1~;~~  C_{k} C_{l} \;=\; - \,C_{l} C_{k}\,,\quad l \neq k\,.
\end{equation}
The action of the symmetric group on this algebra is given by the
natural permutations on the generators\,: $s \cdot C_{k} = C_{s(k)} \cdot s$
for all $s \in S_{n}$.
The algebra $\Sa$ may be regarded \cite{Jones,Stem} as the twisted group
algebra for a certain central extension of the Weyl group of type $B_n$\,.
 
The algebra $\Sa$ has a natural \Z -grading\,: any element
$s \in S_{n}$ has degree zero, while each of the Clifford
generators $C_{1},\ldots,C_{n}$ has degree one. We will use
the notion of \Z -graded irreducibility\,:
a module over the \Z -graded algebra $\Sa$ is
{\it irreducible} if the even part of its supercommutant
equals ${\mathbb C}$.
Furthermore, if the supercommutant coincides with
${\mathbb C}$ then the module is called {\it absolutely irreducible}.

Recall that a partition
$\lambda = (\lambda_{1},\lambda_{2},\ldots,\lambda_{l})$ of $n$
is {\it strict} if all its parts are distinct, that is
$\lambda_{1} > \lambda_{2} > \cdots > \lambda_{l} > 0$\,.
Following the notation introduced by Morris~\cite{Mor},
\mbox{we write $\lambda \stp n$.}
The irreducible modules of the \Z -graded algebra $\Sa$ are parametrised
by the strict partitions $\lambda$ of $n$, while
the irreducible $\Sa$-module labelled by $\lambda$ is
absolutely irreducible if and only if the number $\ell_\lambda$ of parts
in $\lambda$ is even \cite{Serg}.
These two facts can be restated in the classical sense as follows\,:
the irreducible representations over ${\mathbb C}$ of the algebra $\Sa$
are labelled by the pairs
$\left(\lambda,(\pm 1)^{\ell_\lambda}\right)$
where $\lambda \stp n$\qs;
see \cite[Section 1]{Naz} and \cite[Section 9]{Stem}.

\label{sergeev_def}
Let $q$ be an indeterminate over ${\mathbb C}$.
For any positive integer $n$, let $\qSa$ denote the associative algebra with
identity generated over the field ${\mathbb C}(q)$ by elements
$T_{1},T_{2},\ldots,T_{n-1}$ and $C_{1},C_{2},\ldots,C_{n}$ subject to
the following relations.
The generators $T_{1},\ldots,T_{n-1}$ obey the Hecke algebra relations
(\ref{sergeev_def1}),
while the generators $C_{1},\ldots,C_{n}$ satisfy the Clifford algebra
relations (\ref{sergeev_def2}).
Furthermore, there are the cross relations
\begin{equation} \label{sergeev_def3}
\begin{array}{rcl}
T_{k} C_{k} &=& C_{k+1} T_{k}\,; \\
T_{k} C_{k+1} &=& C_{k} T_{k} - (q - q^{-1})(C_{k} - C_{k+1})\,; \\
T_{k} C_{l} &=& C_{l} T_{k}\,, \hspace{4mm} l \neq k,k+1
\end{array}
\end{equation}
for all possible $k,l$.
The algebra $\Sa$ may be recovered as the degenerate case where $q=1$.

We will denote the element $q - q^{-1}$ by $\varepsilon$.
Using the first relation in (\ref{sergeev_def1}), the inverse of the
Hecke generator $T_{k}$ is $T_{k} - \varepsilon$\qs; it is then apparent
that the first two relations in (\ref{sergeev_def3}) are equivalent.
The algebra $\qSa$ has a natural \Z -grading\,: namely, the generators
$T_{1},\ldots,T_{n-1}$ are specified to be even while the Clifford
generators $C_{1},\ldots,C_{n}$ still have degree one.

The Clifford algebra 
has a natural basis
\begin{equation}
\label{natural_basis}
{\cal C} ~=~ \left\{ \left. C_{k_{1}} \cdots \,C_{k_{p}} \right|
1 \leq k_{1} < \cdots < k_{p} \leq n \right\}.
\end{equation}
Given any permutation $s \in S_{n}$\,, let us take any reduced decomposition
$s = s_{j_{1}} \cdots\, s_{j_{r}}$ in terms of the standard Coxeter generators
$s_{1},s_{2},\ldots,s_{n-1}$. Then as usual
define the element $T_{s} \in \qSa$ by
$T_{s} = T_{j_{1}} \cdots\, T_{j_{r}}$.
The second and third relations in (\ref{sergeev_def1}) imply that
this definition does not depend on the reduced decomposition for $s$.

By the defining relations 
(\ref{sergeev_def1}), (\ref{sergeev_def2}) and (\ref{sergeev_def3}),
the elements \,$T_{s} \cdot C \;=\; T_{s} \cdot C_{k_{1}} \cdots \, C_{k_{p}}$
where $s \in S_{n}$ and $C \in {\cal C}$ form a linear basis in the
${\mathbb C}(q)$-algebra $\qSa$.
In particular, we have

\bp \label{dimension}
The algebra\/ $\qSa$ has dimension\/ 
$2^{n} \cdot n\hspace{1pt}!$ over\/ $\CCq$.
\ep

\noindent
The next result is obtained by using the standard technique of \cite{GU}.

\bp \label{semisimple}
The algebra\/ $\qSa$ over\/ $\CCq$ is semisimple.
\ep

\noindent
{\it Proof.}
We will verify that the ${\mathbb C}(q)$-algebra $\qSa$ has a zero radical.
Introduce the associative algebra $\qSast$ generated by the
elements $T_{1},\ldots,T_{n-1}$ and $C_{1},\ldots,C_{n}$ over the ring
${\mathbb C}[\,q,q^{-1}]$ of Laurent polynomials in $q$.
We will view it as an infinite-dimensional algebra over $\mathbb{C}$.
Suppose there exists a non-zero element $R\in\operatorname{rad}\,\qSa$\qs;
we will bring this to a contradiction.
Multiplying $R$ by some element in ${\mathbb C}[q]$, we can assume
that $R\in \qSast$.
On division by a suitable power of $q-1$,
we can assume further that $R\notin(q-1)\cdot \qSast$.

Define the ${\mathbb C}$-algebra homomorphism
$\varpi: \qSast \rightarrow \Sa$ by
$T_{k} \mapsto s_{k},C_{k} \mapsto C_{k},q \mapsto 1$.
We have the equality $s\,C\cdot\varpi(R)=\varpi(T_sCR)$
for any $s\in S_n$ and $C\in\cal C$.
Since $\operatorname{rad}\,\qSa$ is a nilpotent left ideal in $\qSa$
then the element $X\cdot\varpi(R)$ is nilpotent for any $X\in \Sa$.
Thus the left ideal $\Sa\cdot\varpi(R)$ is contained in the radical
of $\Sa$.
But the algebra $\Sa$ is semisimple, hence $\varpi(R) = 0$.
Thus $R\in\operatorname{ker} \varpi = (q-1)\cdot \qSast$.
This is the contradiction.\hspace{3mm}\BOX

\smallskip \noindent
We will prove later that the irreducible $\qSa$-modules are parametrized
by the same strict partitions of $n$ as the irreducible $\Sa$-modules\qs; see
Theorem \ref{V_results}.
We conclude this section by describing the combinatorial objects known
as shifted tableaux\qs; these are analogous to the classical Young tableaux.
The standard reference for these analogues is 
\cite[\mbox{Section III.7}]{Mac}.

Let us fix a strict partition $\lambda = (\lambda_{1},\lambda_{2},\ldots,
\lambda_{l})$ of the integer $n$.
Then the {\it shifted Young diagram of shape} $\lambda$ is the array of
$n$ boxes into $l$ rows with $\lambda_{i}$ boxes in the $i$-th row,
such that each row is shifted by one position to the
right relative to the preceding row.
A {\it shifted tableau of shape} $\lambda$ is an array obtained by
inserting the symbols $1,2,\ldots,n$ bijectively into the $n$ boxes of
the shifted Young diagram for $\lambda$.
We denote an arbitrary shifted tableau by the symbol $\Lambda$.
The set of
all shifted tableaux with shape $\lambda$ is denoted by ${\cal T}_{\lambda}$.
The symmetric group $S_{n}$ acts transitively on the tableaux
$\Lambda \in {\cal T}_{\lambda}$ by permutations on their entries, that is
$$
(s \cdot \Lambda)(i,j) \;=\; s( \Lambda(i,j) ) \hspace{10mm}
\mbox{for all $s \in S_{n} , \Lambda \in {\cal T}_{\lambda}$}.
$$
A shifted tableau is {\it standard} if the symbols increase along
each row (from left to right) and down each column\qs;
the subset of all standard tableaux in ${\cal T}_{\lambda}$
is written as ${\cal S}_{\lambda}$.

There are two distinguished elements in ${\cal S}_{\lambda}$\,:
the {\it row tableau} $\Lambda^{r}$ obtained by inserting the symbols
$1,2,\ldots,n$ consecutively by rows into the shifted diagram of $\lambda$,
and the {\it column tableau} $\Lambda^{c}$ in which the symbols
$1,2,\ldots,n$ occur consecutively by columns.
For example, the row and column tableaux for $\lambda = (4,3,1) \stp 8$ are
given respectively by
\[ \Lambda^{r} = ~\begin{array}{|c|c|c|c|}
\hline 1 & 2 & 3 & 4 \\
\hline \multicolumn{1}{c|}{ } & 5 & 6 & 7 \\
\cline{2-4} \multicolumn{2}{c|}{ } & 8 & \multicolumn{1}{|c}{ } \\
\cline{3-3} \end{array}
\hspace{8mm} \mbox{and} \hspace{8mm}
\Lambda^{c} = ~\begin{array}{|c|c|c|c|}
\hline 1 & 2 & 4 & 7 \\
\hline \multicolumn{1}{c|}{ } & 3 & 5 & 8 \\
\cline{2-4} \multicolumn{2}{c|}{ } & 6 & \multicolumn{1}{|c}{ } \\
\cline{3-3} \end{array}~. \]

Our construction will
involve certain sequences of the integers $1,2,\ldots,n$ derived from these
two special tableaux.
However, the following notation is introduced for an arbitrary standard
tableau $\Lambda \in {\cal S}_{\lambda}$.
Let us denote by $(\Lambda)$ the sequence of integers obtained from the
shifted tableau $\Lambda$ by reading the symbols along the rows from left
to right ordered from the top row to the bottom row.
Similarly, let $(\Lambda)^{\ast}$ denote the sequence derived from $\Lambda$
by reading the symbols down the columns taken from the left column to the
right column.
The special tableaux described above satisfy the obvious properties
$(\Lambda^{r}) = (\Lambda^{c})^{\ast} = ( 1,2,\ldots,n )$.

Now let $\Lambda \in {\cal S}_{\lambda}$ be fixed.
For each index $k=2,\ldots,n$, we denote by ${\cal A}_{k}$ and ${\cal B}_{k}$
the subsequences of $(\Lambda)$ consisting of all entries $j<k$ which occur
respectively after and before $k$ in this sequence.
Similarly, let ${\cal A}_{k}^{\ast}$ and ${\cal B}_{k}^{\ast}$ denote the
subsequences in $(\Lambda)^{\ast}$ consisting of the entries $j<k$ which occur
respectively after and before $k$ in this column sequence.
Denote by $a_{k},b_{k}$ and $a_{k}^{\ast},b_{k}^{\ast}$
the lengths of the sequences ${\cal A}_{k},{\cal B}_{k}$ and
${\cal A}_{k}^{\ast},{\cal B}_{k}^{\ast}$ respectively.

There is a bijection between the set ${\cal T}_{\lambda}$ and the symmetric
group $S_{n}$ described in the following way\,: given any shifted tableau
$\Lambda$ of shape $\lambda$, we define the permutation
$w_{\Lambda} \in S_{n}$~by
$$ 
w_{\Lambda} \;=\; \left( \begin{array}{cccc}
1 & 2 & ~~\cdots~~ & n \\
p_{n} & p_{n-1} & \cdots & p_{1}
\end{array} \right) 
$$
where $(p_{1},\ldots,p_{n})$ is the column
sequence $(\Lambda)^{\ast}$.
This bijection preserves the action of the symmetric group\,:
$w_{s \cdot \Lambda} = s \, w_{\Lambda}$ for all $s \in S_{n}$.
We will \mbox{also require} a second bijection ${\cal T}_{\lambda}
\rightarrow S_{n}$
described as follows\,: given any tableau $\Lambda \in {\cal T}_{\lambda}$
the permutation $s_{\Lambda} \in S_{n}$ is specified by the rule
$s_{\Lambda} \cdot \Lambda = \Lambda^{c}$\,. That is, $s_{\Lambda}$ is the
unique element which maps $\Lambda$ onto the column tableau.
Let $w_0\in S_n$ be the element of maximal length\qs; we have
$w_{0}(k) = n+1-k$ for each $k$.

\bl
\label{perm_product}
The permutations\/ $w_{\Lambda}$ and\/ $s_{\Lambda}$ obey the property\/
$s_{\Lambda} w_{\Lambda} = w_{0}$ for any $\Lambda \in {\cal T}_{\lambda}$.
\el

\noindent
The next result \cite[Lemma 2.4]{Naz} gives reduced decompositions for
$w_{\Lambda}$ and $s_{\Lambda}$ where $\Lambda \in {\cal S}_{\lambda}$.

\bl
\label{perm_decomps}
Given any standard tableau $\Lambda \in {\cal S}_{\lambda}$, there are
reduced decompositions
$$
w_{\Lambda}\,=
\prod_{k=2,\ldots,n}^{\rightarrow}
\left(\,\,
\prod_{p=1,\ldots,b_{k}^{\ast}}^{\rightarrow}s_{k-p}
\right), \hspace{8mm}
s_{\Lambda}\,=
\prod_{k=2,\ldots,n}^{\leftarrow}
\left(\,\,
\prod_{p=1,\ldots,a_{k}^{\ast}}^{\leftarrow}s_{k-p}
\right).
$$
\el

\noindent
Combining these results yields a reduced decomposition for the
element $w_{0}$\,. In particular, for $\Lambda = \Lambda^{c}$
we obtain the obvious reduced decomposition
$$
w_{0} \,=
\prod_{k=2,\ldots,n}^{\rightarrow}
\left(\,\,
\prod_{p=1,\ldots,k-1}^{\rightarrow}s_{k-p}
\right).
$$
The arrow notation on a product indicates the
orientation of its (non-commuting) factors.
\section{The Affine Sergeev Algebra}
\setcounter{equation}{0}
In this section, we introduce a certain object underlying the representation
theory of $\qSa$.
This will be referred to as the {affine Sergeev algebra} and denoted
by $\ASA$\qs; it is a $q$-analogue of the degenerate affine Sergeev algebra
employed in \cite{Naz}.
We will consider representations of $\ASA$ induced from one-dimensional
representations of a certain maximal commutative subalgebra\qs; these are
called {principal series representations}, cf.\ \cite{Rog}.

The {\it affine Sergeev algebra} $\ASA$ is the associative unital algebra
generated over the field ${\mathbb C}(q)$ by $\qSa$ together with the
pairwise-commuting invertible elements $X_{1},X_{2},\ldots,X_{n}$ subject
to the relations
\begin{equation}
\label{affine_rels1}
\begin{array}{rcl}
T_{k} X_{k} &=& X_{k+1} T_{k}-
\varepsilon \left( X_{k+1} - C_{k}C_{k+1} X_{k}\right);
\\
T_{k} X_{k+1} &=& X_{k} T_{k}+
\varepsilon \left( 1 + C_{k}C_{k+1} \right)
X_{k+1}\,;
\\
T_{k} X_{l} &=& X_{l} T_{k}\,,\quad l \neq k,k+1
\end{array}
\end{equation}
for all possible $k,l$ and also the relations
\begin{equation} \label{affine_rels2}
C_{k} X_{k} \;=\; X_{k}^{-1} C_{k}\,;~~ C_{k} X_{l} \;=\; X_{l} C_{k}\,,
\hspace{4mm} l \neq k\,.
\end{equation}
Note that
the first and second relations in (\ref{affine_rels1}) can be deduced
from the single relation
\begin{equation}
\label{Murphy_rel}
\left( T_{k} - \varepsilon \,C_{k}C_{k+1} \right) X_{k} T_{k} \;=\; X_{k+1}.
\end{equation}

The centre $\operatorname{Z}(\ASA)$ will be precisely described
in Proposition \ref{max_comm}.
However, for present purposes, it is sufficient to know that the centre
contains any symmetric polynomial in the elements
$X_1 + X_1^{-1},\ldots,X_n + X_n^{-1}$\qs.
This can be verified by following
\cite[Proposition 3.1]{Naz}.
In particular, the square of the product
\begin{equation}
\label{vander}
\prod_{1 \leq k < l \leq n} 
\left( X_{k} + X_{k}^{-1} - X_{l} - X_{l}^{-1} \right) 
\end{equation}
belongs to the centre $\operatorname{Z}(\ASA)$.
Let $\loc$ denote the localisation of the affine Sergeev algebra by this
central element.
Since each factor in the product (\ref{vander}) can be written as
$$
X_{k}^{-1} \left( X_{k}X_{l} \rule{0mm}{3.5mm} - 1 \right)
\left( X_{k}X_{l}^{-1} - 1 \right),
$$
then it follows that the localisation $\loc$ contains the elements
\begin{equation} \label{local_gen}
\frac{1}{X_{k}X_{l} - 1}~,~~\frac{1}{X_{k}X_{l}^{-1} - 1}~; \qquad
1 \leq k < l \leq n.
\end{equation}
This fact enables us to introduce certain elements in $\loc$ playing a key
role in the sequel.

We will write ${\mathbb C}(q)[X]$ for the subalgebra
${\mathbb C}(q)[X_{1}^{\pm 1},X_{2}^{\pm 1},\ldots,X_{n}^{\pm 1}]$
in $\ASA$\qs;
that is, the space of {Laurent polynomials} in the affine generators
$X_{1},\ldots,X_{n}$.
Let ${\mathbb C}(q)(X)$ denote the subalgebra in the localisation
$\loc$ generated over ${\mathbb C}(q)[X]$ by the elements (\ref{local_gen}).
Now for each index $k=1,2,\ldots,n-1$,
define the element $\Phi_{k} \in \loc$ by
\begin{equation} \label{Phi_def}
\Phi_{k} ~=~ T_{k} + \frac{\varepsilon}{X_{k}X_{k+1}^{-1} - 1} +
\frac{\varepsilon}{X_{k}X_{k+1} - 1} \cdot C_{k}C_{k+1}.
\end{equation}
It follows from the defining relations
(\ref{affine_rels1}), (\ref{affine_rels2})
 that these elements satisfy the properties
\begin{eqnarray}
\Phi_{k} \cdot X_{k} &=& X_{k+1} \cdot \Phi_{k}\,; \nonumber \\
\Phi_{k} \cdot X_{k+1} &=& X_{k} \cdot \Phi_{k}\,; \label{Phi_actions} \\
\Phi_{k} \cdot X_{l} &=& X_{l} \cdot \Phi_{k} ~,~~ l \neq k,k+1\,; \nonumber
\end{eqnarray}
cf.\ \cite{Lus}.
The next result is also established using (\ref{affine_rels1}) and 
(\ref{affine_rels2})\qs; cf.\ \cite[Proposition 3.2]{Naz}.

\bp \label{Phi_rels}
The elements\/ $\Phi_{1},\ldots,\Phi_{n-1}$ obey the following relations
in $\loc$.
\[
\begin{array}{rcl}
\Phi_{k}^{2} &=& \bds 1 - \varepsilon^{2} \cdot \left(
\frac{X_{k}X_{k+1}^{-1}}{(X_{k}X_{k+1}^{-1} - 1)^{2}} +
\frac{X_{k}^{-1}X_{k+1}^{-1}}{(X_{k}^{-1}X_{k+1}^{-1} - 1)^{2}} \right)\eds ;
\\
\Phi_{k} \Phi_{k+1} \Phi_{k} &=& \Phi_{k+1} \Phi_{k} \Phi_{k+1}\, ; \\[1mm]
\Phi_{k} \Phi_{l} &=& \Phi_{l} \Phi_{k}\,, \hspace{4mm} | k-l | > 1.
\end{array}
\]
\ep

\noindent
Given a permutation $s \in S_{n}$ and a reduced decomposition
$s=s_{k_{p}} \cdots \;s_{k_{1}}$
define the element $\Phi_{s} \in \loc$ by
$\Phi_{s} = \Phi_{k_{p}} \cdots \;\Phi_{k_{1}}$.
The second and third relations in Proposition \ref{Phi_rels} imply that this
definition is independent of the reduced decomposition.
The equalities in (\ref{Phi_actions}) demonstrate that the adjoint
action of each element $\Phi_{k} \in \loc$ on
$X_{1},\ldots,X_{n}$ coincides with the standard action of the basic
permutation $s_{k} \in S_{n}$\qs. Thus
\begin{equation}
\label{generalised_action}
\Phi_{s} \cdot X_{k} ~=~ X_{s(k)} \cdot \Phi_{s} ~,\hspace{8mm} k=1,\ldots,n
\end{equation}
for any permutation $s \in S_{n}$.
The family of elements
$\{ \,\Phi_{s} \in \loc \;|\;s \in S_{n}\,\}$
will be used throughout this paper. Let us now prove the following result.

\bp \label{max_comm}
a) The subalgebra\/ $\CS$ in\/ $\ASA$ is maximal commutative.

\noindent
b) The centre\/ $\operatorname{Z}(\ASA)$ of the affine Sergeev
algebra consists precisely of the symmetric polynomials in the elements\/
$X_{k} + X_{k}^{-1}\,,~k=1,\ldots,n$.
\ep

\noindent
{\it Proof.}
a)~We follow an approach by Cherednik \cite[Section 1]{Cher4}.
First, we remark that for any $s\in S_n$ 
the element $\Phi_{s} \in \loc$ has the form
$$
\Phi_{s} \;=~ T_{s} \,+\, \sum_w\,\sum_{C \in {\cal C}}\
f_{w,C} \cdot T_{w}\,C
$$
for some $f_{w,C} \in {\mathbb C}(q)(X)$. Here the first summation is over
the permutations $w \in S_{n}$ with $\mbox{length}(w) < \mbox{length}(s)$.
Thus the elements $\Phi_{s} C$ for all $s \in S_{n}\,,~C \in {\cal C}$ are
linearly independent over ${\mathbb C}(q)(X)$.
In particular, there is a direct sum decomposition
\begin{equation}
\label{local_decomp}
\loc ~= \underset{s \in S_{n},C \in {\cal C}}{\oplus}\
\Phi_{s} C \cdot {\mathbb C}(q)(X) ~=
\underset{s \in S_{n},C \in {\cal C}}{\oplus}\
{\mathbb C}(q)(X) \cdot \Phi_{s} C\,.
\end{equation}

Now consider any element $Z \in \ASA$ in the centralizer of
${\mathbb C}(q)[X]$ as an element in $\loc$. Decompose this element
relative to (\ref{local_decomp}) as
$$
Z \,=\, \sum_{s \in S_{n}}\,\sum_{C \in {\cal C}}\
z_{s,C} \cdot \Phi_{s} C
~;~~ z_{s,C} \in {\mathbb C}(q)(X).
$$
Suppose that $z_{s,C} \neq 0$ for some pair $(s,C) \neq (1,1)$.
Then it follows from (\ref{local_decomp}) that there exists
$P \in {\mathbb C}(q)[X]$ such that $[Z,P] \neq 0$\,: if $C$ is non-trivial
then take $P = X_{j}$ corresponding to any letter $C_{j}$ appearing in the
word $C$. If  $C=1$ and $s$ is non-trivial, then by the
property (\ref{generalised_action}) we take any polynomial in
${\mathbb C}(q)[X]$ not invariant under the action of $s \in S_{n}$.
Thus $Z \in {\mathbb C}(q)(X)$.
However, by definition, $Z\in \ASA$.
Hence $Z \in {\mathbb C}(q)[X]$.

b)~The centre $\operatorname{Z}(\ASA)$ is contained in the maximal
commutative subalgebra ${\mathbb C}(q)[X]$\qs; that is, the central elements
are Laurent polynomials in $X_{1},\ldots,X_{n}$.
We will prove that these polynomials have the stated property.
Firstly, let us check that every symmetric polynomial in the variables
$X_1 + X_1^{-1},\ldots,X_n + X_n^{-1}$
is central in $\ASA$ by showing
that any such polynomial $P$ commutes with each of its generators.

\noindent
i) The polynomial $P$ obviously commutes with the affine generators
$X_{1},\ldots,X_{n}$.\\
ii) The defining relations (\ref{affine_rels2}) show that the element
$X_{k} + X_{k}^{-1}$ commutes with the Clifford generator~$C_{l}$
for any indices $k$ and $l$.\\
iii) Using the definition (\ref{Phi_def}) and (i),(ii) above, the statement
that $P$ commutes with the generators $T_{1},\ldots,T_{n-1}$ is equivalent
to commutation with the elements $\Phi_{1},\ldots,\Phi_{n-1}$.
Since the polynomial $P$ is symmetric, the latter fact follows
from property (\ref{generalised_action}).

Hence, any symmetric polynomial in
$X_1 + X_1^{-1},\ldots,X_n + X_n^{-1}$ is a central element in $\ASA$.
Conversely, the same reasoning as in (a) establishes that the centre
$\operatorname{Z}(\ASA)$ consists precisely of these
elements.\hspace{3mm}\BOX

\smallskip\noindent
Now fix a character $\chi$ of the maximal commutative subalgebra
$\CS\,\subset\,\ASA$ valued in the algebraic closure $\KK$ of the field $\CCq$.
This character is uniquely specified by its values
$\chi(X_{1}),\chi(X_{2}),\ldots,\chi(X_{n})$
where each $\chi(X_{k})$ lies in the multiplicative group
$\KK^\ast$ of $\KK$.
The symmetric group $S_{n}$ has a natural action on the
character $\chi$\,: for any $s \in S_{n}$, we have
$$ s \cdot \chi(X_{k}) \;=\; \chi(X_{s^{-1}(k)})~, \hspace{8mm} k=1,\ldots,n.
$$
Consider the representation $\pi_{\chi}$ over $\KK$
of the algebra $\ASA$ induced from the character $\chi$\,. 
The space $M_{\chi}$ of the representation $\pi_{\chi}$ can be identified
with the finite-dimensional $\KK$-algebra
$\qSdash=\qSa\otimes_{\,\CCq}\KK$.
The generators $T_{1},\ldots,T_{n-1}$ and $C_{1},\ldots,C_{n}$ act in the
representation space $M_{\chi}$ via the usual left multiplication\qs; while the
action of the elements $X_{1},\ldots,X_{n}$ is determined through the defining
relations (\ref{affine_rels1}) and (\ref{affine_rels2}) in $\ASA$ by
$$ X_{k} \cdot m \;=\; (X_{k} m) \cdot 1 \hspace{10mm}
\mbox{for all $m \in M_{\chi}$}\,. $$
Here the action of $X_{k}$ on the identity vector is specified through $\chi$,
that is $X_{k} \cdot 1 = \chi(X_{k})$.

The character $\chi$ is said to be {\it generic} if 
$\chi(X_{k})\neq\chi(X_{l})^{\pm1}$ for all $k\neq l$.
For a generic character $\chi$ the action of $\ASA$ in $M_{\chi}$
extends to each element $\Phi_{s} \in \loc$.
This extended action is also denoted by $\pi_{\chi}$\,.
In the next section, we will see that the elements
$\Phi_{w_{\Lambda}} \in \loc$ with
$\Lambda \in {\cal S}_{\lambda}$
have a well-defined action in $M_{\chi}$ for certain non-generic characters
$\chi$\qs; see Theorem \ref{fusion_theorem}.
For the remainder of this section, we assume that the character
$\chi$ is generic.

\bp
\label{intertwiners}
Given any permutation\/ $s \in S_{n}$ then the operator\/ $\mu_{s}$ of
right multiplication in\/ $\qSdash$ by the element\/ $\pi_{\chi}(\Phi_{s})(1)$
is an intertwining operator $M_{s \cdot \chi}\rightarrow M_{\chi}$.
\ep
\noindent
{\it Proof.}
The action of the generators $T_{1},\ldots,T_{n-1}$ and $C_{1},\ldots,C_{n}$
in the representation spaces $M_{\chi}=M_{s \cdot \chi}=\qSdash$
is through left multiplication and commutes with the operator $\mu_{s}$\,.
It therefore remains for us to verify that the action of the elements
$X_1,\ldots,X_n$ commutes with $\mu_{s}$.
Since the vector $1\in\qSdash$ is cyclic for these actions, then it is
sufficient to demonstrate that the operators
$\pi_{\chi}(X_{k})\cdot\mu_{s}$
and $\mu_{s}\cdot\pi_{s \cdot \chi}(X_{k})$ coincide on the identity vector
for each $k=1,2,\ldots,n$.
This result is established using the property (\ref{generalised_action}) of
the element $\Phi_{s}$\,:
$$
\pi_{\chi}(X_{k}) \left( \pi_{\chi}(\Phi_{s})(1) \right) ~=~
\pi_{\chi}(X_{k} \Phi_{s})(1) ~=~ \pi_{\chi}(\Phi_{s} X_{s^{-1}(k)})(1) ~=~
\pi_{\chi}(\Phi_{s})\left( \pi_{s \cdot \chi}(X_{k})(1) \right)
\hspace{3mm} \BOX
$$

\smallskip \noindent
The proof of Proposition \ref{intertwiners} shows that the element
$\pi_{\chi}(\Phi_{s})(1) \in M_{\chi}$ is an eigenvector for each of the
operators $\pi_{\chi}(X_{k})$\qs;
we will need the following result in Section 6.
\bc
\label{scalar_action}
Given any\/ $s \in S_{n}$ and\/ $k=1,\ldots,n$, we have the equality
$$ \pi_{\chi}(X_{k}) \left( \pi_{\chi}(\Phi_{s})(1) \right) ~=~
(s \cdot \chi)(X_{k}) \cdot \pi_{\chi}(\Phi_{s})(1).
$$
\ec

\noindent
Next, we introduce $q$-analogues for the Jucys-Murphy elements considered
in \cite{Naz}. Using these elements, we will describe a homomorphism
$\iota : \ASA \rightarrow \qSa$.
For each $k=1,2,\ldots,n$, define the {\it Jucys-Murphy element}
$J_{k} \in \qSa$ inductively by
\begin{equation}
\label{Jucys_def}
J_{k} ~=~ \left\{ \begin{array}{ll}
1 & \quad \mbox{for $k=1$\,;} \\
\left( T_{k-1} - \varepsilon \,C_{k-1}C_{k} \right) J_{k-1} T_{k-1} &
\quad \mbox{for $k=2,\ldots,n$.}
\end{array} \right.
\end{equation}
\bp
\label{Murphy_hom}
A homomorphism\/ $\iota : \ASA \rightarrow \qSa$ which is identical
on the subalgebra $\qSa \subset \ASA$ is uniquely specified by
$X_{1} \mapsto 1$.
Then $X_{k} \mapsto J_{k}$ for each $k=1,2,\ldots,n$.
\ep

\noindent
{\it Proof.}
We will verify that the elements $J_1,\ldots,J_n\in\qSa$ satisfy the same
relations with the generators $T_1,\ldots,T_{n-1}$ and $C_1,\ldots,C_n$
as the elements $X_{1},\ldots,X_{n}\in\ASA$ respectively.
By the definition (\ref{Jucys_def})
we have
$T_l\,J_k=J_k\,T_l$
and
$C_l\,J_k=J_k\,C_l$
for any $l>k$\,.
The first two relations in (\ref{affine_rels1}) are equivalent to the
single relation (\ref{Murphy_rel}). But again by the definition
(\ref{Jucys_def}) we have 
$
\left( T_{k} - \varepsilon \,C_{k}C_{k+1} \right) J_{k} T_{k} \;=\; J_{k+1}\,.
$
It now suffices to verify the following relations\,:
\[
\begin{array}{rcl}
T_l\,J_{k+1}=J_{k+1}\,T_l\,,
\quad
l<k\,;
\qquad
C_{k+1}\,J_{k+1}=J_{k+1}^{-1}\,C_{k+1}\,;
\\
C_l\,J_{k+1}=J_{k+1}\,C_l\,,
\quad
J_l\,J_{k+1}=J_{k+1}\,J_l\,,
\qquad
l\leq k\,.
\end{array}
\]
We will verify all these relations
by induction on $k$. The initial case $k=0$ is trivial.

The relation $T_{l} J_{k+1} = J_{k+1} T_{l}$
with $l<k-1$ immediately follows from the inductive assumption.
In the remaining case
$l = k-1$ we have $k \geq 2$ and
$$
J_{k+1}~=~( T_{k} - \varepsilon \,C_{k} C_{k+1} )( T_{k-1} -
\varepsilon \,C_{k-1} C_{k} ) J_{k-1} T_{k-1} T_{k}\,.
$$
The relations (\ref{sergeev_def1}) and
(\ref{sergeev_def2}),(\ref{sergeev_def3}) then provide the equality
\begin{equation}
\label{***} 
T_{k-1} J_{k+1} ~=~ ( T_{k} - \varepsilon \,C_{k} C_{k+1} )
( T_{k-1} - \varepsilon \,C_{k-1} C_{k} ) T_{k} J_{k-1} T_{k-1} T_{k}\,.
\end{equation}
But the element $T_{k}$ commutes with $J_{k-1}$.
The right hand side of (\ref{***}) then equals $J_{k+1} T_{k-1}$ 
by the second relation in (\ref{sergeev_def1}).

By using the definition (\ref{Jucys_def}), the product $C_{k+1}\,J_{k+1}$
can be expressed as $(T_k-\varepsilon)\,C_{k}\,J_{k}\,T_{k}$.
The inductive assumption provides the equality
$C_{k} J_{k} = J_{k}^{-1}C_{k}$. Using the relations (\ref{sergeev_def1})
and (\ref{sergeev_def2}),(\ref{sergeev_def3}) again,
the above expression becomes
$$
(T_{k}-\varepsilon)\,J_{k}^{-1}
\left(T_{k}C_{k+1}+\varepsilon\,(C_{k} - C_{k+1})\right)
\,=\,T_{k}^{-1} J_{k}^{-1}
\left(T_{k}-\varepsilon\,C_{k}C_{k+1}\,\right)^{-1}\,C_{k+1}\,=\,
J_{k+1}^{-1}\,C_{k+1}\,.
$$

The equality $C_l\,J_{k+1}=J_{k+1}\,C_l$ with $l<k$ also follows immediately
from the inductive assumption. Let us check the equality
$C_{k}\,J_{k+1} = J_{k+1}\,C_{k}$.
The product $C_{k}\,J_{k+1}$ can be written as
$( T_{k} - \varepsilon \,C_{k} C_{k+1}) C_{k+1} J_{k} T_{k}$.
But the factors $J_{k}$ and $C_{k+1}$ here commute,
and the result follows by the first relation in (\ref{sergeev_def3}).

Now using (\ref{Jucys_def})
along with the inductive assumption, we get
$J_l\,J_{k+1}=J_{k+1}\,J_{l}$ for each $l<k$.
The remaining case $l=k$ can be settled by writing the element $J_{k}$ as the
product $( T_{k-1} - \varepsilon \,C_{k-1}C_{k} ) J_{k-1} T_{k-1}$.
The element $J_{k+1}$ commutes with each factor here.
\hspace{3mm}\BOX

\smallskip\noindent
Let us now determine the character $\chi$ by introducing an
array with entries in $\KK^\ast$.
Fix an array of shifted shape $\lambda$ with
entries $x(i,j)\in\KK^{\ast}$ and determine the
values $\chi(X_{1}^{-1}),\ldots,\chi(X_{n}^{-1})$ by
$$
w_{\Lambda^{r}} \cdot \chi(X_{k}^{-1}) \,=\, x(i,j)~,
\hspace{5mm} k=\Lambda^{r}(i,j).
$$
In fact, we then have the equality
\begin{equation}
\label{char_def} 
w_{\Lambda} \cdot \chi(X_{k}^{-1}) \,=\, x(i,j)~, \hspace{5mm} k=\Lambda(i,j)
\end{equation}
for any shifted tableau $\Lambda \in {\cal T}_{\lambda}$.
The reason for specifying the values $\chi(X_{k}^{-1})$ rather than the
usual values $\chi(X_{k})$ is purely notational and will become apparent
later on.

Consider the action of the elements
$\Phi_{s} \in \loc$ on the identity vector
$1 \in M_{\chi}$ for any generic character $\chi$.
First, let us present some additional notation.
For any $k=1,2,\ldots,n-1$ let us 
introduce the rational function of $x,y\in\KK$ valued in the algebra $\qSdash$
\begin{equation}
\label{psi_def}
\psi_{k}(x,y) ~=~ T_{k} + \frac{\varepsilon}{x^{-1}y - 1} +
\frac{\varepsilon}{xy - 1} \,C_{k}C_{k+1}\,.
\end{equation}
The action of $\Phi_{1},\ldots,\Phi_{n-1}$ on the identity $1 \in M_{\chi}$ is
determined through the relations (\ref{affine_rels1})
and (\ref{affine_rels2})\,:
for $\chi(X_{k}^{-1}) = x$ and $\chi(X_{k+1}^{-1}) = y$, we have
\begin{equation} \label{psi_intro}
\pi_{\chi}(\Phi_{k})(1) ~=~ \psi_{k}(x,y)\,.
\end{equation}

Let us fix a standard tableau $\Lambda \in {\cal S}_{\lambda}$.
The final result in the present section generalises the equality
(\ref{psi_intro}) by describing how the element
$\Phi_{w_\Lambda}$ acts in the
representation $M_{\chi}$.
It follows directly from Lemma \ref{perm_product} and
Proposition \ref{intertwiners} that
\begin{equation} \label{pi_prod}
\pi_{\chi}(\Phi_{w_{0}})(1) \;=\; \pi_{w_{\Lambda} \cdot \chi}
(\Phi_{s_{\Lambda}})(1) \cdot \pi_{\chi}(\Phi_{w_{\Lambda}})(1).
\end{equation}
Using Lemma \ref{perm_decomps}, we obtain the following decompositions for
the factors on the right hand side in (\ref{pi_prod})\qs;
see \cite[Section 4]{Naz} for a detailed proof of this result.
\bp \label{psi_decomps}
For each\/ $k=1,2,\ldots,n$ let us define\/ $x_{k} = x(i,j)$ where
$k = \Lambda(i,j)$.\/
Then
\begin{eqnarray*}
\pi_{\chi}(\Phi_{w_{\Lambda}})(1) &=& \bds \prod_{k=2,\ldots,n}^{\rightarrow}
\left( \prod_{p=1,\ldots,b_{k}^{\ast}}^{\rightarrow}
\psi_{k-p}(x_{k},x_{{\cal B}_{k}^{\ast}(p)}) \right) , \eds \\
\pi_{w_{\Lambda} \cdot \chi}(\Phi_{s_{\Lambda}})(1) &=& \bds
\prod_{k=2,\ldots,n}^{\leftarrow} \left(
\prod_{q=1,\ldots,a_{k}^{\ast}}^{\leftarrow}
\psi_{k-q}(x_{{\cal A}_{k}^{\ast}(a_{k}^{\ast}-q+1)},x_{k}) \right) . \eds
\end{eqnarray*}
\ep

\noindent
In the next section we study the elements
$\pi_{\chi}(\Phi_{w_{\Lambda}})(1)$
for some non-generic characters $\chi$.


\section{The Elements \boldmath{$\psi_{\Lambda}$} in the
Algebra \boldmath{$\qSprime$} }
\setcounter{equation}{0}
We start this section with introducing certain finite extension $\FF$
of the field $\CCq$. In \text{Section 6} we will show that $\FF$
is a splitting field for the semisimple algebra $\qSa$ over $\CCq$. We write
$$
[m]_{q^2}=\frac{q^{\,2m}-q^{-2m}}{q^{\hskip1pt2}-q^{-2}}
$$
for any integer $m$\qs. Notice that $[0]_{q^2}=0$ and $[1]_{q^2}=1$.
The field $\FF$ is obtained from $\CCq$ by adjoining a square root of
$[m]_{q^2}$ for each $m=2,\ldots,n$.
The $\FF$-algebra $\qSa\otimes_{\,\CCq}\FF$ will be denoted by $\qSprime$,
it is semisimple due to Proposition \ref{semisimple}.
In this section we will define an element $\psi_{\Lambda} \in \qSprime$
for each standard tableau $\Lambda \in {\cal S}_{\lambda}$
as a certain specialization of $\pi_{\chi}(\Phi_{w_{\Lambda}})(1)$.
The element $\psi_{\Lambda^{r}} \in \qSprime$ associated with the row tableau
$\Lambda^{r}$ has particular significance. In Section 6 it will
provide a $q$-analogue of the symmetrizer constructed in \cite{Naz}.

Consider the rational functions $\psi_{k}(x,y)$ of $x,y\in\KK$
defined by (\ref{psi_def}).
These functions are valued in the algebra $\qSdash$.

\bl \label{psi_rel}
The functions\/ $\psi_{k}(x,y)$ satisfy the equations
\[
 \begin{array}{rcll}
\psi_{k}(y,x) \,\psi_{k}(x,y) &=& \bds 1 - \varepsilon^{2}\cdot
\left(
\frac{x^{-1}y}{(x^{-1}y - 1)^{2}} + \frac{xy}{(xy - 1)^{2}}
\right); ~~\eds &
\\
\psi_{k}(x,y) \,\psi_{l}(z,w) &=& \psi_{l}(z,w) \,\psi_{k}(x,y)\,,
\qquad|k-l| \geq 2\,; &
\\[1mm]
\psi_{k}(x,y) \,\psi_{k+1}(z,y) \,\psi_{k}(z,x) &=& \psi_{k+1}(z,x)
\,\psi_{k}(z,y) \,\psi_{k+1}(x,y) &
\end{array}
\]
for all possible\/ $k$ and\/ $l$. Furthermore, we also have the equalities
$$
\psi_{k}(x,y)^{2} \;=\; -\,\varepsilon\cdot\frac{x+y}{x-y}\cdot\psi_{k}(x,y)
\,+\,1\,-\,\varepsilon^{2}\cdot
\left(
\frac{x^{-1}y}{(x^{-1}y - 1)^{2}}+ \frac{xy}{(xy - 1)^{2}}
\right)
$$
and
$$
C_k\,\psi_k(x,y)=\psi_k(x,y^{-1})\,C_{k+1}\,,
\qquad
C_{k+1}\,\psi_k(x,y)=\psi_k(x^{-1},y)\,C_k\,.
$$
\el

\noindent
{\it Proof.}
Let $\chi$ be a generic character of $\CS$ 
such that $\chi(X_{k}^{-1}) = x$ and $\chi(X_{k+1}^{-1}) = y$.
Then it follows from Proposition \ref{intertwiners}
and the equality (\ref{psi_intro}) that
$$ \pi_{\chi}(\Phi_{k}^{2})(1) ~=~ \pi_{s_{k} \cdot \chi}(\Phi_{k})(1) \cdot
\pi_{\chi}(\Phi_{k})(1) ~=~ \psi_{k}(y,x) \,\psi_{k}(x,y). $$
On the other hand, the action of the element $\Phi_{k}^{2}$ on the identity
vector $1 \in M_{\chi}$ may be evaluated explicitly using the first relation in
Proposition \ref{Phi_rels}\,:
$$ \pi_{\chi}(\Phi_{k}^{2})(1) \;=\; 1 - \varepsilon^{2} \cdot \left(
\frac{x^{-1}y}{(x^{-1}y - 1)^{2}} + \frac{xy}{(xy - 1)^{2}} \right) . $$
The second equation in Lemma \ref{psi_rel} is an immediate consequence of
the relations (\ref{sergeev_def1}) to (\ref{sergeev_def3}).
Now let us take a generic character $\chi$ of $\CS$ such that
$\chi(X_{k}^{-1}) = z$,
$\chi(X_{k+1}^{-1}) = x$ and
$\chi(X_{k+2}^{-1}) = y$.
The second relation in Proposition \ref{Phi_rels} gives the equality 
$$ \pi_{\chi}(\Phi_{k} \Phi_{k+1} \Phi_{k})(1) ~=~
\pi_{\chi}(\Phi_{k+1} \Phi_{k} \Phi_{k+1})(1)\,; $$
evaluating the actions on each side using Proposition \ref{intertwiners} gives
\[
 \begin{array}{l}
\pi_{s_{k+1}s_{k} \cdot \chi}(\Phi_{k})(1) \cdot \pi_{s_{k} \cdot \chi}
(\Phi_{k+1})(1) \cdot \pi_{\chi}(\Phi_{k})(1) \\[1.5mm]
\hspace*{11mm} =~ \pi_{s_{k}s_{k+1} \cdot \chi}(\Phi_{k+1})(1) \cdot
\pi_{s_{k+1} \cdot \chi}(\Phi_{k})(1) \cdot \pi_{\chi}(\Phi_{k+1})(1)\,.
\end{array}
 \]
The equality (\ref{psi_intro}) implies that this is precisely the
third equation in Lemma \ref{psi_rel}.
The last three equalities can be established by direct computation\qs;
here the details are omitted.\hspace{3mm}\BOX

\noindent
Note that the third relation in Lemma \ref{psi_rel} is the
{\it Yang-Baxter equation} with the spectral parameters
$x,y$ and $z$\qs; cf.\ \cite{Jim}.
Now consider the following condition on the pair $(x,y)$\,:
\begin{equation}
\label{idem_con}
\frac{x^{-1}y}{(x^{-1}y - 1)^{2}} + \frac{xy}{(xy - 1)^{2}} ~=~
\frac{1}{\varepsilon^{2}}\,\,.
\end{equation}
This constraint is a $q$-analogue to the condition (4.11) in \cite{Naz}\qs;
it will be referred to as the {\it idempotency condition} on $(x,y)$
in view of the following consequence of Lemma~\ref{psi_rel}.
\bc
\label{psi_inv}
a) Suppose that the pair $(x,y)$ satisfies
{\em(\ref{idem_con})} and\/ $y\neq x,x^{-1}$ then
$$
\psi_{k}(x,y)^{2} \;=\; - \varepsilon\cdot\frac{x+y}{x-y}\cdot\psi_{k}(x,y)\,.
$$

\noindent
b) Suppose that the pair $(x,y)$ does not satisfy {\em(\ref{idem_con})}
and\/ $y \neq x,x^{-1}$.
Then\/ $\psi_{k}(x,y)$ is invertible
and
$$ 
\psi_{k}(x,y)^{-1} \,=\; \left[ \,1 - \varepsilon^{2} \cdot \left(
\frac{x^{-1}y}{(x^{-1}y - 1)^{2}} + \frac{xy}{(xy - 1)^{2}} \right) \,\right]
^{-1} \!\psi_{k}(y,x)\,.
$$
\ec

\noindent
Before proceeding any further, let us examine the idempotency condition
(\ref{idem_con}) in more detail.
Consider new variables $(u,v)$ related to $(x,y)$ by
\begin{equation}
\label{scal_sub}
\frac{x+x^{-1}}{2}\;=\;\frac{q\,u^{2}+q^{-1}\,u^{-2}}{q+q^{-1}}\;,
\qquad
\frac{y+y^{-1}}{2}\;=\;\frac{q\,v^{2}+q^{-1}\,v^{-2}}{q+q^{-1}}\,.
\end{equation}
Performing this substitution, the idempotency condition (\ref{idem_con})
takes the form
\begin{equation} \label{quad_con}
\left( q^{2} - u^{2}v^{-2} - u^{-2}v^{2} + q^{-2} \right)
\left( q^{2} - q^{2}u^{2}v^{2} - q^{-2}u^{-2}v^{-2} + q^{-2} \right) \;=~ 0\;;
\end{equation}
this can be factorised as
$$
\left( qu^{2} - q^{-1}v^{2} \right) \left( qu^{-2} - q^{-1}v^{-2} \right)
\left( q^{-1}u^{-2} - q^{-1}v^{2} \right) \left( q^{3}u^{2} - q^{-1}v^{-2}
\right) \;=~ 0\,.
$$
Thus the equation (\ref{quad_con}) has four solutions\,:
$$
v^2\,=\,q^2u^2\,,\quad
v^2\,=\,q^{-2}u^2\,,\quad
v^2\,=\,u^{-2}\,,\quad
v^2\,=\,q^{-4}u^{-2}\,.
$$
Note that the last two solutions can be obtained from the initial two
by the transformation $u \mapsto q^{-1} u^{-1}$.
This transformation exchanges the quadratic factors in (\ref{quad_con})\qs;
whereas the first equality in (\ref{scal_sub}) remains invariant.
Hence each pair $(x,y)$ satisfying the condition (\ref{idem_con}) is
obtained via the substitution (\ref{scal_sub}) from a pair $(u,v)$ obeying
\begin{equation} \label{V_idem_con}
v^2\,=\,q^{\,\pm2}\,u^2.
\end{equation}

Now consider the rational function in the variables $x,y,z$ 
\begin{equation} \label{psi_prod}
\psi_{k}(x,y) \,\psi_{k+1}(z,y) \,\psi_{k}(z,x)
\end{equation}
on the left hand side of the third equation in Lemma \ref{psi_rel}.
This function is regular only when
$y \neq x^{\pm 1} ,~ z \neq y^{\pm 1} ,~ z \neq x^{\pm 1}$.
However, on restriction such that $(x,y)$ satisfies the idempotency
condition, the function (\ref{psi_prod}) has a value at $z=y$.
To formulate this result we introduce some notation.
For $k=1,\ldots,n-2$, define the rational function
\begin{eqnarray*}
\theta_{k}(x,y) &=& \psi_{k}(x,y) \cdot T_{k+1} \cdot \psi_{k}(y,x) ~-~
\varepsilon^{2} \cdot \psi_{k}(x,y) \left[ \frac{x^{-1}y}{(x^{-1}y - 1)^{2}} -
\frac{xy}{(xy - 1)^{2}} \,C_{k}C_{k+1} \right. \\
 & & \hspace{5mm} \left. +\, \frac{1}{(xy - 1)(x^{-1}y - 1)} \,C_{k+1}C_{k+2} +
\frac{1}{(xy - 1)(x^{-1}y - 1)} \,C_{k+2}C_{k} \right]
\end{eqnarray*}
valued in the algebra $\qSdash$.
Let ${\cal I}$ be the subset in $\KK^3$ consisting of all triples $(x,y,z)$
such that the pair $(x,y)$ satisfies (\ref{idem_con}) and $y\neq x,x^{-1}$.
Then we have the following simple lemma.
\bl
\label{sing_rem}
The restriction of the rational function 
{\em(\ref{psi_prod})} to\/ $\cal I$ is regular at\/ $z=y,y^{-1}$.
At $z=y$ this restriction coincides with\/ $\theta_{k}(x,y)$.
\el

\noindent
{\it Proof.}
The function (\ref{psi_prod}) will be written in a form where the singular
component at $z=y,y^{-1}$ can be removed using the constraint (\ref{idem_con}).
The product (\ref{psi_prod}) can be expanded as
\[
\begin{array}{l}
\bds \psi_{k}(x,y) \cdot T_{k+1} \cdot \psi_{k}(z,x) ~+~ \varepsilon \cdot
\psi_{k}(x,y) \cdot T_{k} \left( \frac{1}{z^{-1}y - 1} +
\frac{C_{k}C_{k+2}}{zy - 1} \right) \eds \\[2mm]
\bds \hspace*{5mm} ~+~ \varepsilon^{2} \!\cdot \psi_{k}(x,y) \cdot \left(
\frac{1}{(z^{-1}y - 1)} \frac{1}{(z^{-1}x - 1)} \,+\, \frac{1}{(z^{-1}y - 1)}
\frac{1}{(zx - 1)} \,C_{k}C_{k+1} \right. \eds \\[2mm]
\bds \hspace*{38mm} \left. +\; \frac{1}{(z^{-1}x - 1)} \frac{1}{(zy - 1)}
\,C_{k+1}C_{k+2} \,+\, \frac{1}{(zx - 1)} \frac{1}{(zy - 1)} \,C_{k} C_{k+2}
\right) .\eds
\end{array}
\]
On restriction to ${\cal I}$ the pair $(x,y)$ satisfies (\ref{idem_con})
with $y\neq x,x^{-1}$.
Then $\psi_{k}(x,y) \cdot \psi_{k}(y,x) = 0$ and adding
$$
 - \;\varepsilon\cdot\psi_{k}(x,y)\,\psi_{k}(y,x)\cdot
\left(\frac{1}{z^{-1}y - 1} + \frac{C_{k}C_{k+2}}{zy - 1} \right)
$$
does not alter the value of (\ref{psi_prod}).
Hence the restriction of (\ref{psi_prod}) to ${\cal I}$ coincides with
\[
\begin{array}{l}
\bds \hspace*{-0.5mm}\psi_{k}(x,y) \cdot T_{k+1} \cdot \psi_{k}(z,x) ~-~
\varepsilon^{2} \!\cdot \psi_{k}(x,y) \left[
\frac{x^{-1}z}{(x^{-1}y - 1)(x^{-1}z - 1)} \,-\,
\frac{xz}{(xy - 1)(xz - 1)} \,C_{k}C_{k+1} \right. \eds \hspace{-1mm} \\[3.6mm]
\bds \left. \hspace*{10mm} +\; \frac{1}{(xy - 1)(x^{-1}z - 1)} \,C_{k+1}C_{k+2}
\,+\, \frac{1}{(xz - 1)(x^{-1}y - 1)} \,C_{k+2}C_{k} \right] .\eds
\end{array}
\]
This function is manifestly regular at $z=y,y^{-1}$.
Moreover, at $z=y$ it equals $\theta_{k}(x,y)$.\hspace{3mm}\BOX

\smallskip\noindent
In the previous section, the character $\chi$ has been determined through
(\ref{char_def}) by a shifted array $\{\,x(i,j)\in\KK^\ast\,\}$
of shape $\lambda$.
Let us now fix the standard tableau $\Lambda\in{\cal S}_\lambda$ and, as in
Proposition~\ref{psi_decomps}, write $x_k=x(i,j)$ for $k=\Lambda(i,j)$.
We will denote by $\cal X$ the subset in $(\KK^\ast)^n$ consisting of all
$n$-tuples $(x_1,\ldots,x_n)$ such that $x_k\neq x_l^{\pm 1}$ for $k\neq l$.
These $n$-tuples correspond to the generic characters $\chi$.
Furthermore, denote by $\cal Y$ the subset in $(\KK^\ast)^n$ consisting of all
$(x_1,\ldots,x_n)$ satisfying the following two conditions\,: first, if $k$ and
$l$ occupy different diagonals in $\Lambda$ then $x_k\neq x_l^{\pm 1}$\qs;
secondly, if these two diagonals are not neighbouring
then $(x_k,x_l)$ does not satisfy (\ref{idem_con}). 
Finally, denote by $\cal F$ the subset in $(\KK^\ast)^n$ consisting of all
$(x_1,\ldots,x_n)$ such that for any adjacent entries $k,l$ in the same row
of $\Lambda$ then the pair $(x_k,x_l)$ satisfies the idempotency condition
(\ref{idem_con}) with $x_k+x_l\neq0$.

If $k=\Lambda(i,j)$
the difference $j-i$ is called the {\it content} of the box in
the shifted diagram $\lambda$ occupied by the symbol $k$ in $\Lambda$.
Let us fix the {\it special\/} point $\qseq \in(\KK^\ast)^n$ where
\begin{equation}
\label{q_def}
q_{\qs k}\,=\, 
[\hskip1pt j\!-\!i\!+\!1]_{q^2}-[\hskip1pt j\!-\!i\hskip1pt]_{q^2}-
\varepsilon\,
\sqrt{\,
[\hskip1pt j\!-\!i\!+\!1]_{q^2}[\hskip1pt j\!-\!i\hskip1pt]_{q^2}
}\,\,,
\qquad
k=\Lambda(i,j)\,.
\end{equation}
Observe that $\qseq \notin {\cal X}$,
while $\qseq \in {\cal F}\cap{\cal Y}$.
Here the membership of ${\cal F}$ follows directly from the
definition (\ref{q_def})\qs; see  (\ref{scal_sub}) and (\ref{V_idem_con}).
Meanwhile, for arbitrary entries $k = \Lambda(i,j)$ and $l = \Lambda(i',j')$,
the definition (\ref{q_def}) implies that 
$$
\left( q_{\qs k}^{-1} q_{\qs l} - 1 \right)
\left( \rule{0mm}{3.5mm}q_{\qs k} - q_{\qs l}^{-1} \right)
\;=\;
\frac{2}{q+q^{-1}}\,\cdot 
\left(
q^{2(j'-i')+1} - q^{2(j-i)+1}+
q^{-2(j'-i')-1} - q^{-2(j-i)-1}
\right)
$$
where $j-i$ and $j'-i'$ are non-negative integers.
This equality shows that $q_{\qs l} = q_{\qs k}$ or $q_{\qs l}
= q_{\qs k}^{-1}$
if and only if $j-i = j'-i'$, that is $k$ and $l$ occupy the same
diagonal in $\Lambda$.
Thus the special point $\qseq$ lies within ${\cal Y}$.

Now for each $\Lambda \in {\cal S}_{\lambda}$ let us define
$$
\psi_{\Lambda}(x_1,\ldots,x_n) \;=\; \pi_{\chi}(\Phi_{w_{\Lambda}})(1)
$$
for any generic character $\chi$.
We consider $\psi_\Lambda(x_1,\ldots,x_n)$ as a rational function
of the variables $x_1,\ldots,x_n$ valued in the algebra $\qSdash$.
By Proposition \ref{psi_decomps}, we have
\begin{equation}
\label{psi_Lam}
\psi_{\Lambda}(x_1,\ldots,x_n) ~=\,
\prod_{k=2,\ldots,n}^{\rightarrow} \left(
\prod_{p=1,\ldots,b_{k}^{\ast}}^{\rightarrow}
\psi_{k-p}(x_{k},x_{{\cal B}_{k}^{\ast}(p)})
\right)
\end{equation}
where the sequences ${\cal B}_{k}^{\ast}$ are as defined in Section 2.
This rational function may have poles outside the set $\cal X$.
However, we will establish that its restriction to $\cal F$ is regular in
$\cal F\cap\cal Y$.
The continuation of the function (\ref{psi_Lam}) to the point $\qseq$ along
$\cal F$ is called the {\it fusion procedure,\/}
this notion has been introduced by Cherednik \cite{Cher2}.

\bt \label{fusion_theorem}
For any standard tableau $\Lambda \in {\cal S}_{\lambda}$\,,
the restriction of\/ $\psi_{\Lambda}\xseq$
to $\cal F$ is regular in\/ $\cal F\cap\cal Y$.
This restriction does not vanish at the point\/ $\qseq$.
\et

\noindent
{\it Proof.}
We will follow the constructive proof from \cite[Theorem 5.6]{Naz}.
Firstly, we remark that Lemma \ref{perm_decomps} gives
the reduced decomposition
$$ 
w_{\Lambda} \,= \prod_{k=2,\ldots,n}^{\rightarrow} \left(
\prod_{p=1,\ldots,b_{k}^{\ast}}^{\rightarrow} s_{k-p} \,\right)
$$
for the element $w_{\Lambda} \in S_{n}$.
By expanding the product (\ref{psi_Lam}) using the definition of the
functions
$\psi_{k}(x,y)$, we obtain a sum with leading term 
$$ 
T_{w_{\Lambda}} = \prod_{k=2,\ldots,n}^{\rightarrow} \left(
\prod_{p=1,\ldots,b_{k}^{\ast}}^{\rightarrow} T_{k-p} \,\right) 
$$
while the remaining terms involve elements $T_{s}$ with
permutations $s \in S_{n}$ of smaller length.
Therefore if the restriction of
$\psi_{\Lambda}\xseq$ is regular at a point, its value must be non-zero.

Our consideration may be restricted to the case for $\Lambda=\Lambda^{c}$.
Namely, using the
equalities in (\ref{pi_prod}) and Proposition \ref{psi_decomps}, we have
\begin{equation}
\label{eqnA}
\prod_{k=2,\ldots,n}^{\leftarrow}
\left(
\prod_{q=1,\ldots,a_{k}^{\ast}}^{\leftarrow}
\psi_{k-q}(x_{{\cal A}_{k}^{\ast}(a_{k}^{\ast}-q+1)},x_{k})
\right)
\cdot\psi_{\Lambda}(x_1,\ldots,x_n)
~=~
\psi_{\Lambda^{c}}
(x_1^{\hskip1pt\prime},\ldots\hskip-.3pt,x_n^{\hskip1pt\prime}\,)
\end{equation}
where we write $x_k^{\hskip1pt\prime}=x(i,j)$ if $k=\Lambda^c(i,j)$.
Each factor $\psi_{k-q}(x_{{\cal A}_{k}^{\ast}(a_{k}^{\ast}-q+1)},x_{k})$
in the product
on the left hand side of (\ref{eqnA}) is regular in $\cal Y$ and has
invertible values.
Hence we can assume that $\Lambda = \Lambda^{c}$.
Then $x_k^{\hskip1pt\prime}$ coincides with $x_k$ for each $k$.

The function $\psi_{\Lambda^{c}}(x_1,\ldots,x_n)$ can be written
as the product
\begin{equation}
\label{eqnB}
\psi_{\Lambda^{c}}(x_1,\ldots,x_n)\;=\;
\theta_{\Lambda^{c}}(x_1,\ldots,x_n)\cdot
\theta^{\,\prime}_{\Lambda^{c}}(x_1,\ldots,x_n)
\end{equation}
where we denote
\begin{eqnarray}
\theta_{\Lambda^{c}}(x_1,\ldots,x_n) &=&
\prod_{k=2,\ldots,n}^{\rightarrow} \left(
\prod_{p=1,\ldots,b_{k}}^{\rightarrow} \psi_{k-p}(x_{k},x_{{\cal B}_{k}(p)})
\right),
\label{eqnC} \\
\theta^{\,\prime}_{\Lambda^{c}}(x_1,\ldots,x_n) &=&
\prod_{k=2,\ldots,n}^{\leftarrow}
\left(
\prod_{q=1,\ldots,a_{k}}^{\leftarrow}
\psi_{n-k+q}(x_{k},x_{{\cal A}_{k}(a_{k}-q+1)})
\right).
\label{eqnD}
\end{eqnarray}
Here, the sequences ${\cal B}_{k}$ and ${\cal A}_{k}$ are defined
for $\Lambda^{c}$ as described in Section 2.
Each factor in the product (\ref{eqnD}) is regular in $\cal Y$.
The proof is complete once we verify that the restriction of the function
$\theta_{\Lambda^{c}}(x_1,\ldots,x_n)$ to $\cal F$
is regular in $\cal F\cap\cal Y$.
The product (\ref{eqnC}) may contain factors
which become singular within $\cal Y\setminus\cal X$.
Specifically, these are the factors corresponding to pairs $(k,p)$
where the symbols $k$ and ${\cal B}_{k}(p)$ stand on the same diagonal
in the tableau $\Lambda^{c}$.
Let us call any such pair {\it singular\,}.
Given any singular pair $(k,p)$, we note that the symbols ${\cal B}_{k}(p)$
and ${\cal B}_{k}(p+1)$ are adjacent within some row of $\Lambda^c$.

Let $\hat{\theta}_{\Lambda^{c}}(x_1,\ldots,x_n)$ denote the product obtained
from (\ref{eqnC}) by inserting the expression
\begin{equation}
\label{ins_exp}
\varepsilon^{-1} \cdot
\frac{x_{{\cal B}_{k}(p)}-x_{{\cal B}_{k}(p+1)}}
{x_{{\cal B}_{k}(p)}+x_{{\cal B}_{k}(p+1)}}
\cdot \psi_{k-p-1}(x_{{\cal B}_{k}(p+1)},x_{{\cal B}_{k}(p)})
\end{equation}
before each factor $\psi_{k-p}(x_{k},x_{{\cal B}_{k}(p)})$ with
singular $(k,p)$.
On restriction to $\xseq \in {\cal F}$, the pair
$(x_{{\cal B}_{k}(p+1)},x_{{\cal B}_{k}(p)})$ satisfies the
condition (\ref{idem_con}) with
$x_{{\cal B}_{k}(p+1)} \neq x_{{\cal B}_{k}(p)}^{\pm 1}$\,.
Then the product
$$
\psi_{k-p-1}(x_{{\cal B}_{k}(p+1)},x_{{\cal B}_{k}(p)})
\,\psi_{k-p}(x_{k},x_{{\cal B}_{k}(p)})
\,\psi_{k-p-1}(x_{k},x_{{\cal B}_{k}(p+1)})
$$
is regular at $x_k=x_{{\cal B}_{k}(p)}$ by Lemma \ref{sing_rem}.
Using Lemma \ref{psi_rel}, it can be shown that this procedure
of inserting the normalised factors (\ref{ins_exp})
into the product (\ref{eqnC})
does not alter its restriction to $\cal F$. 
On the other hand, the rational function
$\hat{\theta}_{\Lambda^{c}}(x_1,\ldots,x_n)$
is regular in $\cal F\cap\cal Y$.\hspace{3mm}\BOX

\smallskip \noindent
Observe that for any tableau $\Lambda\in{\cal S}_{\lambda}$
the coordinates $q_1,\ldots,q_n$ belong to the subfield $\FF\subset\KK$.
Hence for each $\Lambda\in{\cal S}_{\lambda}$ we can define 
the element $\psi_{\Lambda} \in \qSprime$ as the value at $\qseq$
of the restriction to $\cal F$ of $\psi_{\Lambda}(x_1,\ldots,x_n)$.
Note that 
the proof of Theorem \ref{fusion_theorem} together with Lemma \ref{sing_rem}
provides an explicit formula for the element $\psi_{\Lambda}$\qs.

Let us make another important observation.
Fix an index $k\in\{1,\ldots,n\}$. Consider the point
$\qseq^\prime\in\FF$ obtained from 
$\qseq$ by inverting the coordinate $q_{\qs k}$\,.
We have $\qseq^\prime\in\cal F\cap\cal Y$.
Again using Theorem \ref{fusion_theorem} define the element
$\psi_\Lambda^{\,\prime}\in\qSprime$  
as the value at $\qseq^\prime$
of the restriction to $\cal F$ of $\psi_{\Lambda}(x_1,\ldots,x_n)$.
By using the last two equalities in Lemma \ref{psi_rel}
along with the definition (\ref{psi_Lam}) we obtain the following proposition.
 
\bp
\label{c_inter}
We have the equality\/ 
$C_k\,\psi_\Lambda=\psi_\Lambda^{\,\prime}\,C_{w_\Lambda^{-1}(k)}$
in the algebra $\qSprime$.
\ep

\noindent
If 
$$
x\,=\,[a+1]_{q^2}-[a]_{q^2}-\varepsilon\,\sqrt{\,[a+1]_{q^2}[a]_{q^2}}\,\,,
\quad
y\,=\,[b+1]_{q^2}-[b]_{q^2}-\varepsilon\,\sqrt{\,[b+1]_{q^2}[b]_{q^2}}
\hskip-15pt
$$
for some non-negative integers $a\neq b$ then at $q\to1$
the element $\psi_k(x,y)\in\qSprime$
degenerates to
$$
s_{k}+
\Bigl(\sqrt{a(a+1)}-\sqrt{b(b+1)}\,\,\Bigr)^{-1}\!\!-
\Bigl(\sqrt{a(a+1)}+\sqrt{b(b+1)}\,\,\Bigr)^{-1}
C_k C_{k+1}
\in\Sa\!\!
$$
by definition (\ref{psi_def}).
Degenerations at $q\to1$ of the elements $\psi_\Lambda\in\qSprime$
were studied in \cite{Naz}. In the subsequent two sections we will give
$q$-analogues of these results from \cite{Naz}. 


\section{On the Divisibility of the Element \boldmath{$\psi_{\Lambda^{r}}$}}
\setcounter{equation}{0}
The opening part of this section examines the left-divisibility of the element
$\psi_{\Lambda^{r}} \in \qSprime$ by certain elements in the algebra $\qSprime$
corresponding to pairs of adjacent row entries in $\Lambda^{r}$.
In the next proposition, the scalars $q_{\qs k}$ and $q_{\qs k+1}$
are defined by (\ref{q_def}) for $\Lambda = \Lambda^{r}$.

\bp
\label{row_ann}
Suppose that\/ $k = \Lambda^{r}(i,j)$ and\/ $k+1 = \Lambda^{r}(i,j+1)$ are
adjacent entries in some row of\/ $\Lambda^{r}$.
Then\/ $\psi_{k}(q_{\qs k},q_{\qs k+1}) \cdot \psi_{\Lambda^{r}} \,=\; 0$.
\ep

\noindent
{\it Proof.}
For any $n$-tuple $\xseq$ in ${\cal X} \cap {\cal F}$
we have
$
\psi_{\Lambda^{r}} \xseq ~=~ \pi_{\chi}(\Phi_{w_{\Lambda^{r}}})(1)
$
where $\chi$ is the generic character determined by $x_{1},\ldots,x_{n}$
through (\ref{char_def}) with $\Lambda = \Lambda^{r}$.
Now the entry $k$ precedes $k+1$ in the column sequence $(\Lambda^{r})^{\ast}$.
By definition of the permutation $w_{\Lambda^{r}} \in S_{n}$\,,
it follows that
$\mbox{length}(w_{\Lambda^{r}}) \;=\;
\mbox{length}(s_{k} w_{\Lambda^{r}}) + 1$.
Then Proposition \ref{intertwiners} gives
\begin{eqnarray}
\pi_{\chi}(\Phi_{w_{\Lambda^{r}}})(1) &=&
\pi_{s_{k}w_{\Lambda^{r}} \cdot \chi}(\Phi_{k})(1) \cdot
\pi_{\chi}(\Phi_{s_{k}w_{\Lambda^{r}}})(1) \nonumber \\
 &=& \psi_{k}(x_{k+1},x_{k}) \cdot \pi_{\chi}(\Phi_{s_{k}w_{\Lambda^{r}}})(1)
\label{row_div1}
\end{eqnarray}
where the second equality is given by (\ref{psi_intro}).
Furthermore, by the definition of the set ${\cal F}$, the pair
$(x_{k+1},x_{k})$ satisfies the idempotency condition (\ref{idem_con})\qs;
thus the first relation in Lemma \ref{psi_rel} gives the
equality $\psi_{k}(x_{k},x_{k+1}) \cdot \psi_{k}(x_{k+1},x_{k}) \;=\; 0$.
It follows from (\ref{row_div1}) that
\begin{equation}
\label{psi_ann}
\psi_{k}(x_{k},x_{k+1}) \cdot \psi_{\Lambda^{r}} \xseq ~=~ 0\;,
\qquad \xseq \in {\cal X} \cap {\cal F}.
\end{equation}
Since the restriction of the rational function on the left hand side of
(\ref{psi_ann}) to ${\cal F}$ is regular in ${\cal F} \cap {\cal Y}$
then the stated result follows by continuation along ${\cal F}$
to the point $\qseq$.\hspace{3mm}\BOX

\smallskip\noindent
An immediate consequence of (\ref{psi_def}) is
\begin{equation}
\label{minus}
\psi_k(y,x) \;=\; \psi_k(x,y) \,+\, \varepsilon\cdot\frac{x+y}{x-y}\,.
\end{equation}
Using this identity we obtain the following corollary to Proposition
\ref{row_ann}.

\bc \label{row_div}
Suppose that\/ $k = \Lambda^{r}(i,j)$ and\/ $k+1 = \Lambda^{r}(i,j+1)$ are
adjacent entries in some row of the tableau $\Lambda^{r}$.
Then
$$
\psi_{k}(q_{\qs k+1},q_{\qs k}) \cdot \psi_{\Lambda^{r}} \;=~ \varepsilon \cdot
\frac{q_{\qs k}+q_{\qs k+1}}{q_{\qs k}-q_{\qs k+1}} \cdot \psi_{\Lambda^{r}}\,.
$$
\ec

\smallskip\noindent
An analogue of Proposition \ref{row_ann} exists for adjacent column
entries in the tableau $\Lambda^{r}$, but it is not as apparent\,:
given any $n$-tuple $\xseq$ in ${\cal F}$, the pair
$(x_{k},x_{l})$ where $k$ and $l$ are adjacent entries in the same column of
$\Lambda^{r}$ does not necessarily satisfy the idempotency condition.
The remainder of this section is devoted to proving this analogue
(Corollary \ref{col_div3}).

Once again, consider the restriction of the rational function
$\psi_{k}(x,y) \,\psi_{k+1}(z,y) \,\psi_{k}(z,x)$
to the set ${\cal I}$.
Lemma \ref{sing_rem} demonstrates that this restriction is regular at
$z=y$ and has identical values to the function $\theta_{k}(x,y)$.
This expression is not necessarily divisible on the right by $\psi_{k}(y,x)$\,.
However, let us consider the rational function
\begin{equation}
\label{right_div}
\psi_{k}(x,y) \,\psi_{k+1}(z,y) \,\psi_{k}(z,x) \cdot \psi_{k}(x,z)
\end{equation}
and denote by $d(x,y)$ the rational function
$$
\varepsilon^{3} \!\cdot y \left\{
(y^{2}-1)\left( \frac{x^{3}}{(xy-1)^{4}}+\frac{x^{-3}}{(x^{-1}y-1)^{4}} \right)
\,+\,\frac{x^{3}-x}{(xy-1)^{4}}
\,+\,\frac{x^{-3}-x^{-1}}{(x^{-1}y-1)^{4}}\,\right\}
$$
valued in $\KK$. We will need the following lemma.

\bl \label{theta_result}
The restriction of the function {\em(\ref{right_div})} to\/ $\cal I$
coincides at\/ $z=y$ with the restriction of the product\/
$d(x,y)\,\psi_{k}(x,y)$.
\el

\noindent
{\it Proof.}
First, let us examine the product $\psi_{k}(z,x) \,\psi_{k}(x,z)$.
The first relation in Lemma \ref{psi_rel} and the idempotency condition
on $(x,y)$ give the equality
$$
\psi_{k}(z,x) \,\psi_{k}(x,z) ~=~ \varepsilon^{2} \cdot \left(
\frac{xy}{(xy - 1)^{2}} + \frac{x^{-1}y}{(x^{-1}y - 1)^{2}} -
\frac{xz}{(xz - 1)^{2}} - \frac{x^{-1}z}{(x^{-1}z - 1)^{2}} \right)
$$
on $\cal I$.
A direct calculation shows that the restriction of (\ref{right_div})
to ${\cal I}$ can be written as
$$
\psi_{k}(x,y) \left( T_{k+1} + \frac{\varepsilon}{z^{-1}y - 1} +
\frac{\varepsilon}{zy - 1}\,C_{k+1}C_{k+2} \right) \cdot \varepsilon^{2} \cdot
\left\{ (yz - 1)\rule{0mm}{3.8mm}(y^{-1}z - 1) \right\} y\cdot\, d(x,y,z)
$$
where $d(x,y,z)$ denotes the sum
\begin{eqnarray*}
& & 
\frac{x^{3}}{(xy - 1)^{2}(xz - 1)^{2}}+
\frac{x^{3} - x}{(xy - 1)^{2}(xz - 1)^{2}(yz - 1)}+
\\
& &
\frac{x^{-3}}{(x^{-1}y - 1)^{2}(x^{-1}z - 1)^{2}}+
\frac{x^{-3} - x^{-1}}{(x^{-1}y - 1)^{2}(x^{-1}z - 1)^{2}(yz - 1)}\,.
\end{eqnarray*}
Thus the restriction to ${\cal I}$ of (\ref{right_div}) is regular at $z=y$
provided that $y^{2} \neq 1$\qs; moreover, it takes identical values to the
function
\[
\begin{array}{l}
\varepsilon^{3} \!\cdot d(x,y,y) ~y (y^{2}-1) \cdot
\psi_{k}(x,y) ~=\\[2mm]
\varepsilon^{3} \!\cdot y \bds\left\{ (y^{2}-1) \left(
\frac{x^{3}}{(xy-1)^{4}} + \frac{x^{-3}}{(x^{-1}y-1)^{4}} \right) +
\frac{x^{3}-x}{(xy-1)^{4}} + \frac{x^{-3}-x^{-1}}{(x^{-1}y-1)^{4}}
\right\} \eds \cdot \psi_{k}(x,y)\,.
\end{array}
\]
As a consequence, the restriction is also regular at $y^{2}=1$.\hspace{3mm}\BOX

\bc \label{theta_zero}
The restriction of\/ {\em (\ref{right_div})} to\/ ${\cal I}$
vanishes at\/ $y=z=1$.
\ec

\noindent
{\it Proof.}
On the set ${\cal I}$, the function $d(x,y)$ takes the value zero
when $y=1$.\hspace{3mm}\BOX

\medskip \noindent
We will use Corollary \ref{theta_zero} with $y = x_{k}$ and $z = x_{l}$
where both $k$ and $l$ stand on the
{\it leading diagonal} of the tableau $\Lambda$\,:
this is the diagonal with the entries $\Lambda(i,i)$.
At the special point $\qseq$ the definition (\ref{q_def}) then gives
$q_{\qs k}=q_{\qs l}=1$.
Next, let us consider the involutive antiautomorphism
$\alpha : \qSprime \rightarrow \qSprime$ defined by
$T_{k}\mapsto T_{n-k}$ and $C_{k}\mapsto C_{n-k+1}$\,.
Note that
$\alpha \left( \psi_{k}(x,y)\rule{0mm}{3.5mm} \right) \,=\, \psi_{n-k}(x,y)$
for each index $k=1,\ldots,n-1$.
The next result concerns the element $\psi_{\Lambda^{c}}$ associated 
with the column tableau.

\bp \label{fixedpoint}
The element\/ $\psi_{\Lambda^{c}} \in \qSprime$ is invariant under the
antiautomorphism $\alpha$.
\ep
\noindent
{\it Proof.}
For the column tableau $\Lambda^{c}$, each sequence ${\cal B}_{k}^{\ast}$
is the complete interval $(1,2,\ldots,k-1)$.
Hence the equality of rational functions (\ref{psi_Lam}) becomes
$$
\psi_{\Lambda^{c}}\xseq ~= \prod_{k=2,\ldots,n}^{\rightarrow} \left(
\prod_{j=1,\ldots,k-1}^{\rightarrow} \psi_{k-j}(x_{k},x_{j}) \right) .
$$
A direct application of the second and third equalities in Lemma \ref{psi_rel}
establishes that
$$
\psi_{\Lambda^{c}}\xseq ~=\; \prod_{k=2,\ldots,n}^{\leftarrow} \left(
\prod_{j=1,\ldots,k-1}^{\leftarrow} \psi_{n-k+j}(x_{k},x_{j}) \right) \;=~
\alpha \left( \rule{0mm}{3.4mm}\, \psi_{\Lambda^{c}}\xseq \,\right) \,.
$$
In particular, the restriction of the function $\psi_{\Lambda^{c}}\xseq$
onto ${\cal X} \cap {\cal F}$ is invariant under $\alpha$.
Now Proposition \ref{fixedpoint} follows by continuation
along ${\cal F}$ to the point $\qseq$.\hspace{3mm}\BOX

\smallskip \noindent
In the proof of Theorem \ref{fusion_theorem} we showed that
the restriction of the function $\theta_{\Lambda^{c}}\xseq$ to ${\cal F}$
is regular at the point $\qseq$.
Define the element $\theta_{\Lambda^{c}} \in \qSprime$ as the value of
the restriction at this point.
We now present the main result in this section\qs; here $q_{\qs k}$ and
$q_{\qs k+1}$ are defined by (\ref{q_def}) for $\Lambda = \Lambda^{c}$.

\bt \label{col_div1}
Suppose that\/ $k = \Lambda^{c}(i,j)$ and\/ $k+1 = \Lambda^{c}(i+1,j)$
are adjacent entries in some column of the column tableau.
Then\/ $\theta_{\Lambda^{c}} \in \qSprime$ is divisible on the left by\/
$\psi_{k}(q_{\qs k+1},q_{\qs k})$.
\et

\noindent
{\it Proof.}
First, observe that Theorem \ref{col_div1} is implied by its
particular case where $k = n-1$. Indeed,
given an arbitrary index $k$, let $\Omega^{c}$
be the tableau obtained from $\Lambda^{c}$ by removing each of the symbols
$k+2,\ldots,n$\,. Then the tableau $\Omega^{c}$ is the column tableau for a
certain partition $\omega \stp k+1$ and
$\theta_{\Lambda^{c}} =\, \theta_{\Omega^{c}} \cdot \tilde{\theta}$
for some element $\tilde{\theta} \in \qSprime$.

Thus we will assume that $k=n-1$.
Let us demonstrate left divisibility of the element $\theta_{\Lambda^{c}}$ by
$\psi_{n-1}(q_{\qs n},q_{\qs n-1})$.
In particular, we will establish the equality
\begin{equation}
\label{theta_div}
\psi_{n-1}(q_{\qs n},q_{\qs n-1}) \cdot \theta_{\Lambda^{c}} \,=~
\varepsilon \cdot \frac{q_{\qs n-1} + q_{\qs n}}{q_{\qs n-1} - q_{\qs n}}
\cdot \theta_{\Lambda^{c}}\,.
\end{equation}
By (\ref{minus}), this is equivalent to verifying
\begin{equation}
\label{theta_ann}
\psi_{n-1}(q_{\qs n-1},q_{\qs n}) \cdot \theta_{\Lambda^{c}} \,=~ 0\,.
\end{equation}

Recall that $n-1 = \Lambda^{c}(i,j)$ and $n = \Lambda^{c}(i+1,j)$.
Let $p_{1} < p_{2} < \cdots < p_{j-i}=n$ be the entries
in the $(i+1)$-th row of the tableau $\Lambda^{c}$.
The restriction to ${\cal F}$ of the factor
$\theta^{\,\prime}_{\Lambda^{c}}\xseq$
in (\ref{eqnB}) is regular at the point $\qseq$\qs;
let us denote its value by $\theta'_{\Lambda^{c}} \in \qSprime$.
Using the second and third equations in Lemma \ref{psi_rel}, this
element can be shown to satisfy
$$
\theta'_{\Lambda^{c}} \cdot \psi_{1}(q_{\qs n-1},q_{\qs n}) ~=
\prod_{k=1,\ldots,j-i}^{\rightarrow}
\psi_{j-i+n-b_{n}-k}(q_{\qs n-1},q_{\qs p_{k}}) ~\cdot ~\eta
$$
for some element $\eta \in \qSprime$\qs; cf.\ \cite[Proposition 2.8\/]{Naz2}.
Then
$$
\psi_{\Lambda^{c}} \cdot \psi_{1}(q_{\qs n-1},q_{\qs n}) ~=~
\theta_{\Lambda^{c}} \cdot \prod_{k=1,\ldots,j-i}^{\rightarrow}
\psi_{j-i+n-b_{n}-k}(q_{\qs n-1},q_{\qs p_{k}}) ~\cdot ~\eta
$$
and it will be sufficient to verify that
\begin{equation}
\label{theta_ann2}
\theta_{\Lambda^{c}} \cdot \prod_{k=1,\ldots,j-i}^{\rightarrow}
\psi_{j-i+n-b_{n}-k}(q_{\qs n-1},q_{\qs p_{k}}) ~=~ 0\,.
\end{equation}
Applying the antiautomorphism $\alpha$ to the equality
\mbox{$\psi_{\Lambda^{c}} \cdot \psi_{1}(q_{\qs n-1},q_{\qs n}) = 0$}
using Proposition~\ref{fixedpoint}, we obtain
$\psi_{n-1}(q_{\qs n-1},q_{\qs n}) \cdot \psi_{\Lambda^{c}} =\, 0$.
Since the element $\psi_{\Lambda^{c}}$ can be realised by multiplying
$\theta_{\Lambda^{c}}$ on the right by the invertible element
$\theta'_{\Lambda^{c}}$, this last equality is equivalent to
the required statement (\ref{theta_ann}).
It therefore remains to establish (\ref{theta_ann2}).
This will be proved using induction on the integer $j-i$.

I. Suppose that $j-i = 1$\qs; that is, the symbol $n$ stands on
the leading diagonal in the column tableau $\Lambda^{c}$. Then $b_{n} = n-1$.
Let $m = \Lambda^{c}(i,i)$ then ${\cal B}_{n}(n-2) = m$
and ${\cal B}_{n}(n-1) = n-1$.
The expression (\ref{eqnC}) defining the function
$\theta_{\Lambda^{c}}\xseq$ takes the form of the product
$$
\theta\xseq \cdot \psi_{2}(x_{n},x_{m})\, \psi_{1}(x_{n},x_{n-1})
$$
where
$$
\theta\xseq=\!\!
\prod_{k=2,\ldots,n-1}^{\rightarrow}
\left(\,
\prod_{p=1,\ldots,b_{k}}^{\rightarrow}\!\psi_{k-p}(x_{k},x_{{\cal B}_{k}(p)})
\right)
\times\!\!
\prod_{p=1,\ldots,n-3}^{\rightarrow}\!\psi_{n-p}(x_{n},x_{{\cal B}_{n}(p)})
$$
on restriction to ${\cal F}$ is regular at $\qseq$.
Due to Lemma \ref{psi_rel} and Corollary \ref{psi_inv}(a),
this restriction satisfies the equality
$$
\theta\xseq \;=\; \varepsilon^{-1} \cdot
\frac{x_{m} - x_{n-1}}{x_{m} + x_{n-1}} \cdot
\theta\xseq \cdot \psi_{1}(x_{n-1},x_{m})\,.
$$
Thus the left hand side in (\ref{theta_ann2}) is the value at $\qseq$
taken by the restriction of
$$
\varepsilon^{-1} \cdot \frac{x_{m} - x_{n-1}}{x_{m} + x_{n-1}} \cdot
\theta\xseq \times
\psi_{1}(x_{n-1},x_{m})\, \psi_{2}(x_{n},x_{m})\, \psi_{1}(x_{n},x_{n-1})
\cdot \psi_{1}(x_{n-1},x_{n})\,.
$$
For $\xseq \in {\cal F}$, the pair $(x_{n-1},x_{m})$ satisfies
the condition (\ref{idem_con}) with \mbox{$x_{m} + x_{n-1} \neq 0$}\qs;
while $q_{\qs m} = q_{\qs n} = 1$.
Corollary \ref{theta_zero} demonstrates that the restriction
vanishes at $\qseq$.

II. Next, consider the case $j-i > 1$.
We will demonstrate that the restriction to ${\cal F}$ of
\begin{equation}
\label{theta_expr}
\theta_{\Lambda^{c}}\xseq \;\cdot \prod_{k=1,\ldots,j-i}^{\rightarrow}
\psi_{j-i+n-b_{n}-k}(x_{n-1},x_{p_{k}})
\end{equation}
has the value zero at $\qseq$.
In this instance, let us specify $m = \Lambda^{c}(i+1,j-1)$ then
$m-1 = \Lambda^{c}(i,j-1)$.
Let $\Omega^{c}$ denote the tableau obtained from $\Lambda^{c}$ by removing
each of the entries $m+1,\ldots,n$\qs; evidently, $\Omega^{c}$ is the column
tableau for a certain partition $\omega\stp m$.
Now the expression in (\ref{eqnC}) defining the function
$\theta_{\Lambda^{c}}\xseq$ can be expanded as
\[
\begin{array}{l}
\theta_{\Omega^{c}}\xseq \cdot \theta'\xseq \,
\psi_{n-b_{n-1}-1}(x_{n-1},x_{m-1}) ~~\times 
\\[2mm]
\hspace*{6mm}\theta''\xseq \, \psi_{n-b_{n}+j-i}(x_{n},x_{m-1})\,
\psi_{n-b_{n}+j-i-1}(x_{n},x_{n-1})\; \theta\xseq
\end{array}
\]
where the functions $\theta'\xseq$ and $\theta''\xseq$ are defined by
\[
\begin{array}{c}
\bds
\theta'\xseq\,=\!\!
\prod_{k=m+1,\ldots,n-2}^{\rightarrow} \!\!
\left(
\prod_{p=1,\ldots,b_{k}}^{\rightarrow}\!
\psi_{k-p}(x_{k},x_{{\cal B}_{k}(p)})
\!\right)
\times \hspace{-2mm}\prod_{p=1,\ldots,b_{n-1}-1}^{\rightarrow}
\hspace{-3.5mm}\!\psi_{n-1-p}(x_{n-1},x_{{\cal B}_{n-1}(p)}),
\eds
\medskip
\\
\bds
\theta''\xseq\,=\,
\prod_{p=1,\ldots,b_{n}-j+i-1}^{\rightarrow}
\psi_{n-p}(x_{n},x_{{\cal B}_{n}(p)})
\eds
\end{array}
\]
while
$$
\theta\xseq ~=
\prod_{k=1,\ldots,j-i-1}^{\rightarrow}
\psi_{n-b_{n}+j-i-1-k}(x_{n},x_{p_{k}})\,.
$$

Since the integers $b_{n-1}$ and $b_{n}$ are related by the equality
$b_{n} = b_{n-1}+j-i$, we have the inequalities
$n - p \,\geq\, n - b_{n-1} + 1$ for each $p = 1,\ldots,b_{n}-j+i-1$\,.
Thus the function $\psi_{n-b_{n-1}-1}(x_{n-1},x_{m-1})$ commutes with
$\theta''\xseq$ by the second relation in Lemma \ref{psi_rel}.
Therefore the function (\ref{theta_expr}) can be written as
\[
\begin{array}{l}
\theta_{\Omega^{c}}\xseq \cdot \theta'\xseq\, \theta''\xseq ~~\times
\\[1.5mm]
\hspace*{5mm} \psi_{n-b_{n-1}-1}(x_{n-1},x_{m-1})\,
\psi_{n-b_{n-1}}(x_{n},x_{m-1})\, \psi_{n-b_{n-1}-1}(x_{n},x_{n-1})
\\[1mm]
\hspace*{5mm} \bds \times~~ \theta\xseq \,\cdot
\prod_{k=1,\ldots,j-i}^{\rightarrow}
\psi_{n-b_{n-1}-k}(x_{n-1},x_{p_{k}}) \eds
\end{array}
\]
while further use of the second and third equalities in
Lemma \ref{psi_rel} gives the expression
\begin{equation}
\label{theta_expr2}
\begin{array}{l}
\theta_{\Omega^{c}}\xseq \cdot \theta'\xseq\, \theta''\xseq ~~\times
\\[1.5mm]
\hspace*{5mm}\bds \psi_{n-b_{n-1}-1}(x_{n-1},x_{m-1})
\,\psi_{n-b_{n-1}}(x_{n},x_{m-1}) \,\psi_{n-b_{n-1}-1}(x_{n},x_{n-1})
\,\psi_{n-b_{n-1}-1}(x_{n-1},x_{n}) \eds \hspace{-4mm}
\\[1mm]
\hspace*{5mm}\bds \times ~\prod_{k=1,\ldots,j-i-1}^{\rightarrow}
\psi_{n-b_{n-1}-1-k}(x_{n-1},x_{p_{k}}) \prod_{k=1,\ldots,j-i-1}^{\rightarrow}
\psi_{n-b_{n-1}-k}(x_{n},x_{p_{k}}) \,. \eds
\end{array}
\end{equation}
It follows from the proof of Theorem \ref{fusion_theorem} that
the product in the first line of (\ref{theta_expr2})
on restriction to ${\cal F}$ is regular at $\qseq$.
The restriction of each of the factors in the
last line is also regular at this point.
Furthermore, for any $\xseq \in {\cal F}$, the pair $(x_{n-1},x_{m-1})$
satisfies the condition (\ref{idem_con}).
Thus Lemma \ref{theta_result} demonstrates that the restriction of the
function on the second line in (\ref{theta_expr2}) coincides at $\qseq$
with the restriction to ${\cal F}$ of the function
$$
d(x_{n-1},x_{m-1}) \cdot \psi_{n-b_{n-1}-1}(x_{n-1},x_{m-1}).
$$
Hence the restriction of (\ref{theta_expr2}) to ${\cal F}$ has the same
value at $\qseq$ as the restriction of
\[
\begin{array}{cl}
 & \bds \theta_{\Omega^{c}}\xseq \cdot \theta'\xseq\, \theta''\xseq \cdot
d(x_{n-1},x_{n}) \,\psi_{n-b_{n-1}-1}(x_{n-1},x_{n}) ~~\times \eds
\\[1mm]
 & \bds \hspace{5mm}\prod_{k=1,\ldots,j-i-1}^{\rightarrow}
\psi_{n-b_{n-1}-1-k}(x_{n-1},x_{p_{k}}) \prod_{k=1,\ldots,j-i-1}^{\rightarrow}
\psi_{n-b_{n-1}-k}(x_{n},x_{p_{k}}) \eds
\\[1.5mm]
=& \bds \theta_{\Omega^{c}}\xseq \cdot \theta'\xseq\, \theta''\xseq \cdot
\prod_{k=1,\ldots,j-i-1}^{\rightarrow} \psi_{n-b_{n-1}-1-k}(x_{n},x_{p_{k}})
~~\times \eds
\\[-1mm]
 & \bds \hspace{5mm}\prod_{k=1,\ldots,j-i-1}^{\rightarrow}
\psi_{n-b_{n-1}-k}(x_{n-1},x_{p_{k}}) \cdot d(x_{n-1},x_{n})\,
\psi_{n-b_{n}}(x_{n-1},x_{n})\eds \,.
\end{array}
\]
Replacing the variable $x_{n}$ by $x_{m-1}$ in every factor
$\psi_{n-b_{n-1}-1-k}(x_{n},x_{p_{k}})$ within the fourth component
of the latter expression does not affect the value of the restriction
at $\qseq$.
Let us denote this modified component by
$$
\bar{\theta}\xseq ~= \prod_{k=1,\ldots,j-i-1}^{\rightarrow}
\psi_{n-b_{n-1}-1-k}(x_{m-1},x_{p_{k}})\,.
$$
For each index $k=1,2,\ldots,j-i-1$, we have the inequalities
$(n-b_{n-1}-1-k) + 2 \,<\, n-p$ for every $1 \leq p \leq b_{n}-j+i-1$\qs;
thus $\bar{\theta}\xseq$ commutes with $\theta''\xseq$.

The proof is complete once we demonstrate that the restriction to ${\cal F}$ of
\begin{equation} \label{theta_expr3}
\theta_{\Omega^{c}}\xseq\, \theta'\xseq\, \bar{\theta}\xseq
\end{equation}
vanishes at the point $\qseq$.
The factors in the function $\theta'\xseq$ are arranged with respect to
the ordering specified by the sequences ${\cal B}_{k}$.
Let us now rearrange these factors in the following manner\,:
for each index $k > m$ appearing in the $(j-1)$-th column of
$\Lambda^{c}$, change the subsequence
$m-1\,,\,p_{\qs 1}\,,\,p_{\qs 2}\,, \ldots ,\,p_{\qs j-i-1}$ in ${\cal B}_{k}$
to
$p_{\qs 1}\,,\,p_{\qs 2}\,, \ldots ,\,p_{\qs j-i-1}\,,\,m-1$.
Recall that
$$
m-1=\Lambda^{c}(i,j-1)\,,
\quad
p_{\qs 1}=\Lambda^{c}(i+1,i+1)\,,
\quad\ldots\,,\quad
p_{\qs j-i-1}=\Lambda^{c}(i+1,j-1).
$$
Let $\theta^{\ast}\xseq$ denote the product
obtained from $\theta'\xseq$ by this rearrangement; put
$$
\tilde{\theta}\xseq ~= \prod_{k=1,\ldots,j-i-1}^{\rightarrow}
\psi_{m-b_{n-1}+i-k}(x_{m-1},x_{p_{k}})\,.
$$
The equality
$$
\theta'\xseq \cdot \bar{\theta}\xseq
\;=\;
\tilde{\theta}\xseq \cdot \theta^{\ast}\xseq
$$
of rational functions
can be established
by using the second and third relations in Lemma \ref{psi_rel}.
Hence the product (\ref{theta_expr3}) becomes
$\theta_{\Omega^{c}}\xseq \,\tilde{\theta}\xseq \,\theta^{\ast}\xseq$.
Furthermore, any singular factors in the restriction of the function
$\theta^{\ast}\xseq$ to ${\cal F}$
can be dealt with as described in the proof of
Theorem \ref{fusion_theorem}\qs; that is, we consider an expression
$\hat{\theta}^{\ast}\xseq$ obtained by inserting normalised factors
into $\theta^{\ast}\xseq$ at certain positions.
Then
$$
\theta_{\Omega^{c}}\xseq \,\tilde{\theta}\xseq \,\theta^{\ast}\xseq ~=~
\theta_{\Omega^{c}}\xseq \,\tilde{\theta}\xseq \,\hat{\theta}^{\ast}\xseq
$$
on ${\cal F}$, while the restriction of $\hat{\theta}^{\ast}\xseq$
to ${\cal F}$ is regular at $\qseq$.
Meanwhile, the inductive assumption on the tableau $\Omega^{c}$ shows that the
restriction to ${\cal F}$ of the function \vspace{-1mm}
$$
\theta_{\Omega^{c}}\xseq \cdot \tilde{\theta}\xseq ~=~
\theta_{\Omega^{c}}\xseq \;\cdot
\prod_{k=1,\ldots,j-i-1\vspace{-1mm}}^{\rightarrow}
\psi_{m-b_{m}+j-1-i-k}(x_{m-1},x_{p_{k}})
$$
vanishes at this special point\qs; thus the restriction of (\ref{theta_expr3})
also vanishes.\hspace{3mm}\BOX

\smallskip \noindent
The next result is a consequence of Proposition \ref{fixedpoint} and
Theorem \ref{col_div1}.
Here $q_{\qs k}$ and $q_{\qs k+1}$ are again
defined by (\ref{q_def}) for the tableau $\Lambda=\Lambda^{c}$.

\bc
\label{col_div2}
Suppose that\/ $k = \Lambda^{c}(i,j)$ and\/ $k+1 = \Lambda^{c}(i+1,j)$ are
adjacent entries within the same column of\/ $\Lambda^{c}$.
Then the element\/ $\psi_{\Lambda^{c}} \in \qSprime$ is\,:

a) divisible on the left by\/ $\psi_{k}(q_{\qs k+1},q_{\qs k})$ and

b) divisible on the right by\/ $\psi_{n-k}(q_{\qs k+1},q_{\qs k})$.
\ec

\noindent
Furthermore, we obtain the following analogue to Corollary \ref{row_div}.
In this instance, we specify $q_{\qs k}$ and $q_{\qs l}$ 
by (\ref{q_def}) for the row tableau $\Lambda=\Lambda^{r}$.

\bc \label{col_div3}
Suppose that\/ $k = \Lambda^{r}(i,j)$ and\/ $l = \Lambda^{r}(i+1,j)$ are
adjacent entries within the same column of the row tableau.
Then
$$
\psi_{\Lambda^{r}} \cdot \psi_{n-p}(q_{\qs l},q_{\qs k}) ~=\;
\varepsilon \cdot \frac{q_{\qs k} + q_{\qs l}}{q_{\qs k} - q_{\qs l}}
\cdot \psi_{\Lambda^{r}}
$$
where\/ $p = \Lambda^{c}(i,j)$ is the entry in the column tableau occupying
the same position as\/ $k$ in $\Lambda^{r}$.
\ec

\noindent
{\it Proof.}
Since $p = \Lambda^{c}(i,j)$ then $p+1=\Lambda^{c}(i+1,j)$.
Thus Corollary \ref{col_div2}(b) yields the equality
$$
\psi_{\Lambda^{c}} \cdot
\psi_{n-p}(q_{\qs p+1}^{\qs \prime},q_{\qs p}^{\qs \prime}) ~=\;
\varepsilon \cdot
\frac{q_{\qs p}^{\qs \prime} + q_{\qs p+1}^{\qs \prime}}
{q_{\qs p}^{\qs \prime} - q_{\qs p+1}^{\qs \prime}}
\cdot \psi_{\Lambda^{c}}
$$
where $q_{\qs p}^{\qs \prime}$ and $q_{\qs p+1}^{\qs \prime}$
are defined by (\ref{q_def}) for $\Lambda=\Lambda^{c}$.
The identification
$
q_{\qs p}^{\qs \prime} = q_{\qs k}\,,~q_{\qs p+1}^{\qs \prime} = q_{\qs l}
$
is obvious. Meanwhile, since
$\psi_{\Lambda^{r}} = \,\mu \cdot \psi_{\Lambda^{c}}$ for some
invertible element $\mu \in \qSprime$ then $\psi_{\Lambda^{c}}$ may be
replaced by $\psi_{\Lambda^{r}}$ in the above equality to give the
stated result.\hspace{3mm}\BOX

\smallskip \noindent
The element $\psi_{\Lambda^{r}}T_{w_{\Lambda^{r}}}^{\,-1}\!\in\qSprime$ is
our analogue of the product $E_\omega\in H_n(q)$ as described in Section~1.
We will justify this claim further in the \mbox{next section.}
We also conjecture that this element of $\qSprime$
is an idempotent up to a multiplier from $\FF^\ast$.
Using Propositions \ref{commutant} and \ref{gamma_inter}
one can prove that the degeneration at $q\to1$ of this element
of $\qSprime$ is an idempotent
in the algebra $\Sa$ up to a certain multiple from
${\mathbb C}^{\hskip1pt\ast}$.
This confirms \text{\cite[Conjecture 9.4]{Naz}.}


\section{The \boldmath{$\qSprime$}-Module \boldmath{$V_{\lambda}$}}
\setcounter{equation}{0}
In this final section for each partition $\lambda\stp n$
we will construct a certain $\qSprime$-module $V_{\lambda}$\qs.
The $\qSdash$-module $V_\lambda\otimes_{\,\FF}\KK$ will appear 
as a subrepresentation in the principal series
representation $M_{\chi_{0}}$ of $\ASA$.
Here $\chi_{0}$ is the character of $\CS$ determined
through any shifted tableau $\Lambda \in {\cal T}_{\lambda}$ by
$
w_{\Lambda} \cdot \chi_{0}(X_{k}^{-1})=q_{\qs k}\,,
$
see the definition (\ref{q_def}).
The subalgebra $\qSa\subset\ASA$ acts in
$M_{\chi_{0}}=\qSdash$ via left multiplication.
We will observe that the action of $\ASA$ in $V_{\lambda}\otimes_{\,\FF}\KK$
factors through the homomorphism $\iota : \ASA \rightarrow \qSa$
introduced in Proposition \ref{Murphy_hom}.
The $\qSprime$-module $V_{\lambda}$ is reducible\hskip1pt; its irreducible
components are described in Theorem \ref{V_results}.

Consider the partial {\it Bruhat order\/} $\succ$ on the elements
in the symmetric group $S_{n}$\,: given permutations $s,s' \in S_{n}$ then
$s \succ s'$ if and only if there are adjacent transpositions
$s_{k_{1}},\ldots,s_{k_{p}}$ such that
$s = s_{k_{p}} \cdots \, s_{k_{1}} \cdot s'$ where
$\mbox{length}(s) = \mbox{length}(s') + p$\,.
For any $\Lambda \in {\cal S}_{\lambda}$ the elements $w_{\Lambda}$ and
$s_{k} w_{\Lambda}$ in $S_{n}$ are neighbours with respect to this
partial ordering for each $k = 1,\ldots,n-1$.
We now examine further the elements $\psi_{\Lambda} \in \qSprime$
introduced in Section 4.

\bp \label{psi_formulae}
Let\/ $\Lambda \in {\cal S}_{\lambda}$ and\/ $k \in \{ 1,2,\ldots,n-1 \}$.
Define\/ $q_k,q_{k+1}\!\in\FF$ by {\em(\ref{q_def})}.
Then
\begin{enumerate}
\smallsep
\item if\/ $s_{k} \cdot \Lambda \in {\cal S}_{\lambda}$ and\/
$s_{k}w_{\Lambda} \succ w_{\Lambda}$ then ~$\psi_{k}(q_{\qs k},q_{\qs k+1})
\cdot \psi_{\Lambda} \,=\, \psi_{s_{k} \cdot \Lambda}$\,;

\item if\/ $s_{k} \cdot \Lambda \in {\cal S}_{\lambda}$ and\/
$s_{k}w_{\Lambda} \prec w_{\Lambda}$ then ~$\psi_{k}(q_{\qs k+1},q_{\qs k})
\cdot \psi_{s_{k} \cdot \Lambda} \,=\, \psi_{\Lambda}$\,;

\item if\/ $s_{k} \cdot \Lambda \not\in {\cal S}_{\lambda}$ then
~$\psi_{k}(q_{\qs k},q_{\qs k+1}) \cdot \psi_{\Lambda} =\, 0$\,.
\end{enumerate}
\ep

\noindent
{\it Proof.}
Let us fix the standard tableau $\Lambda \in {\cal S}_{\lambda}$ and the
index $k$.
The tableau \mbox{$s_{k} \cdot \Lambda$} obtained by interchanging the entries
$k$ and $k+1$ within $\Lambda$, may or may not be standard.
Suppose that $k = \Lambda(i,j)$ and $k+1 = \Lambda(i',j')$.
Since the tableau $\Lambda$ is standard, we have four possibilities
to consider\,:
\begin{quote}
$i' > i,~j' < j$\,;
\hspace{7.5mm}
$i' < i,~j' > j$\,;
\hspace{7.5mm}
$i' = i,~j' = j+1$\,;
\hspace{7.5mm}
$i' = i+1,~j' = j$.
\end{quote}
The tableau \mbox{$s_{k} \cdot \Lambda$} is standard in
precisely the first two cases.
Furthermore, using the definition of the permutation $w_{\Lambda} \in S_{n}$
gives $w_{s_{k} \cdot \Lambda} \succ w_{\Lambda}$ in the first situation while
$w_{s_{k} \cdot \Lambda} \prec w_{\Lambda}$ in the second instance.
We will now examine each case separately.

i)~In the first case, we have $s_{k} \cdot \Lambda \in {\cal S}_{\lambda}$
and $s_{k}w_{\Lambda} \succ w_{\Lambda}$.
Let us verify the equality in (a) for the elements $\psi_{\Lambda}$
and $\psi_{s_{k} \cdot \Lambda}$.
Given any $n$-tuple $\xseq \in {\cal X} \cap {\cal F}$, we have
$$
\psi_{\Lambda}\xseq \;=\; \pi_{\chi}(\Phi_{w_{\Lambda}})(1)
$$
where
$\chi$ is the generic character determined through (\ref{char_def})
by $x_{1},\ldots,x_{n}$.
Similarly for the tableau $s_{k} \cdot \Lambda$, we have the equality
$$
\psi_{s_{k} \cdot \Lambda}\xseq \;=\; \pi_{\chi}(\Phi_{s_{k}w_{\Lambda}})(1)
$$
on ${\cal X} \cap {\cal F}$.
Using Proposition \ref{intertwiners}, the restriction of
$\psi_{s_{k} \cdot \Lambda}\xseq$ to ${\cal F}$ can be written as
$$ \pi_{w_{\Lambda} \cdot \chi}(\Phi_{k})(1)
\cdot \pi_{\chi}(\Phi_{w_{\Lambda}})(1) ~=~ \psi_{k}(x_{k},x_{k+1}) \cdot
\psi_{\Lambda}\xseq . $$
Here the equality follows from (\ref{psi_intro}) since
$w_{\Lambda} \cdot \chi (X_{k}^{-1}) = x_{k}\,,~w_{\Lambda} \cdot
\chi (X_{k+1}^{-1}) = x_{k+1}$.
The stated result is obtained by continuation along ${\cal F}$ to the
point $\qseq$.

ii)~The conditions in (b) are realised in the second case\,:
$s_{k} \cdot \Lambda \in {\cal S}_{\lambda}$ and
$s_{k}w_{\Lambda} \prec w_{\Lambda}$.
Let us consider the tableau $\Lambda' = s_{k} \cdot \Lambda$.
Then $s_{k} \cdot \Lambda' = \Lambda$ while
$s_{k}w_{\Lambda'} \succ w_{\Lambda'}$.
The result follows from (a) for the tableau $\Lambda'$.

iii)~In the third instance, the entries $k$ and $k+1$ are adjacent in some
row of $\Lambda$ and hence $s_{k} \cdot \Lambda \not\in {\cal S}_{\lambda}$.
The equality in (c) has been verified for the row tableau
$\Lambda^{r}$ in Proposition \ref{row_ann}\qs; this proof
extends directly to an arbitrary tableau $\Lambda \in {\cal S}_{\lambda}$.

iv)~In the final case, the equality in (c) has been established for
the column tableau $\Lambda^{c}$ in the preceding section.
Let us now extend this equality to an arbitrary tableau
$\Lambda \in {\cal S}_{\lambda}$.
The entries $k = \Lambda(i,j)$ and $k+1 = \Lambda(i+1,j)$ are adjacent
within some column of $\Lambda$.
Since the symbol $k$ precedes $k+1$ in the column sequence $(\Lambda)^{\ast}$
then $s_{k}w_{\Lambda} \prec w_{\Lambda}$
by the definition of the element $w_{\Lambda} \in S_{n}$\qs;
it also follows from Lemma \ref{perm_product} that
$s_{\Lambda}s_{k} \succ s_{\Lambda}$.
Furthermore, the definition of the permutation $s_{\Lambda}$ gives
$s_{\Lambda}s_{k} = s_{p}s_{\Lambda}$ where $p = \Lambda^{c}(i,j)$ and
$p+1 = \Lambda^{c}(i+1,j)$ are the entries in $\Lambda^{c}$ occupying the
same positions as $k$ and $k+1$ respectively, in the tableau $\Lambda$.
In particular, we have
$\mbox{length}(s_{p}s_{\Lambda}) \,=\, \mbox{length}(s_{\Lambda}s_{k}) \,=\,
\mbox{length}(s_{\Lambda}) + 1$.

For any generic character $\chi$, the element
$\pi_{w_{\Lambda} \cdot \chi}(\Phi_{s_{\Lambda}s_{k}})(1) \in \qSprime$
can be expressed as
$$
\pi_{s_{k}w_{\Lambda} \cdot \chi}(\Phi_{s_{\Lambda}})(1) \cdot
\pi_{w_{\Lambda} \cdot \chi}(\Phi_{k})(1) ~=~
\pi_{s_{k}w_{\Lambda} \cdot \chi}(\Phi_{s_{\Lambda}})(1) \cdot
\psi_{k}(x_{k},x_{k+1})
$$
using Proposition \ref{intertwiners} and the equality (\ref{psi_intro}).
Note that $x_{k}$ and $x_{k+1}$ are determined by $\chi$ via (\ref{char_def})
for the tableau $\Lambda$.
The same element
$\pi_{w_{\Lambda} \cdot \chi}(\Phi_{s_{\Lambda}s_{k}})(1)$
may also be expressed as
\begin{eqnarray*}
\pi_{s_{\Lambda}w_{\Lambda} \cdot \chi}(\Phi_{p})(1) \cdot
\pi_{w_{\Lambda} \cdot \chi}(\Phi_{s_{\Lambda}})(1) &=&
\pi_{w_{0} \cdot \chi}(\Phi_{p})(1) \cdot
\pi_{w_{\Lambda} \cdot \chi}(\Phi_{s_{\Lambda}})(1) \\
 &=& \psi_{p}(x_{k},x_{k+1}) \cdot
\pi_{w_{\Lambda} \cdot \chi}(\Phi_{s_{\Lambda}})(1)
\end{eqnarray*}
where the second equality is given by (\ref{psi_intro}) since
$w_{0} \cdot \chi (X_{p}^{-1}) = x_{k}$ and
$w_{0} \cdot \chi (X_{p+1}^{-1}) = x_{k+1}$.
Thus we have verified the identity
$$
\psi_{p}(x_{k},x_{k+1}) \cdot
\pi_{w_{\Lambda} \cdot \chi}(\Phi_{s_{\Lambda}})(1) ~=~
\pi_{s_{k}w_{\Lambda} \cdot \chi}(\Phi_{s_{\Lambda}})(1) \cdot
\psi_{k}(x_{k},x_{k+1})
$$
for any generic $\chi$.
Meanwhile the equality (\ref{pi_prod}) gives
$$
\psi_{\Lambda^{c}}\xseq ~=~
\pi_{w_{\Lambda} \cdot \chi}(\Phi_{s_{\Lambda}})(1) \cdot \psi_{\Lambda}\xseq
$$
for any $n$-tuple $\xseq \in {\cal X}$.
Combining these two equalities on ${\cal X}$, we obtain
\begin{equation}
\label{iden2}
\psi_{p}(x_{k},x_{k+1}) \cdot \psi_{\Lambda^{c}}\xseq ~=~
\pi_{s_{k}w_{\Lambda} \cdot \chi}(\Phi_{s_{\Lambda}})(1) \cdot
\psi_{k}(x_{k},x_{k+1}) \cdot \psi_{\Lambda}\xseq\,.
\end{equation}
In the proof of Theorem \ref{col_div1}, we have established that
the restriction to ${\cal F}$ of the rational function on the left
hand side in (\ref{iden2}) vanishes at the point $\qseq$.
Furthermore, it can be shown that the factor
$\pi_{s_{k}w_{\Lambda} \cdot \chi}(\Phi_{s_{\Lambda}})(1)$ on the right hand
side in (\ref{iden2}) is regular in ${\cal Y}$ and has invertible values\qs;
cf.\ \cite[Proposition 7.1]{Naz}.
Thus the equality in (c) follows by the continuation of (\ref{iden2})
along ${\cal F}$ to the special point $\qseq$.\hspace{3mm}\BOX

\smallskip\noindent
Now let $V_{\lambda} = \qSprime\,\psi_{\Lambda^{r}}$ be the left ideal in
the algebra $\qSprime$ generated by the element $\psi_{\Lambda^{r}}$.
Consider $\qSdash$ as the space of the principal series
representation $M_{\chi_{0}}$ of the affine Sergeev algebra $\ASA$.
Theorem \ref{V_theorem} will demonstrate that this action of the affine
algebra in $\qSdash$ preserves the subspace $V_\lambda\otimes_{\,\FF}\KK$.
Before stating this result, let us introduce some additional notation.
The Clifford algebra, now viewed as the subalgebra in $\qSprime$ generated
over ${\FF}$ by the elements $C_{1},\ldots,C_{n}$ has the natural
basis ${\cal C}$ described in (\ref{natural_basis}).
For each element $C = C_{k_{1}} \cdots \,C_{k_{p}} \in {\cal C}$
let us define the function
$\nu_{C} : \{ 1,2,\ldots,n \} \rightarrow \{\,1,-1\,\}$ by
$$ 
\nu_{C}(k) ~=~ \left\{
\begin{array}{cl}
+1 & \quad\mbox{if the word $C$ contains the letter $C_{k}$\,;} \\
-1 & \quad\mbox{otherwise.}
\end{array} \right.
$$
Proposition \ref{Murphy_hom} describes a homomorphism
$\iota : \ASA \rightarrow \qSa$ which is identical on the
subalgebra $\qSa$.
The part (c) of the following theorem shows that the action of the
affine algebra $\ASA$ in the subspace
$V_\lambda\otimes_{\,\FF}\KK\subset M_{\chi_{0}}$
factors through this homomorphism.

\bt
\label{V_theorem}
a) The elements\/ $C \psi_{\Lambda}$ where\/ $C \in {\cal C}$ and\/
$\Lambda \in {\cal S}_{\lambda}$ form a\/ ${\FF}$-basis in $V_{\lambda}$.

\noindent
b) For each\/ $k = 1,2,\ldots,n$, the action of the element\/ $X_{k}$ on the
basis vectors is given by
$$ \pi_{\chi_{0}}(X_{k})(C \psi_{\Lambda}) ~=~ q_{\qs k}^{\qs \nu_{C}(k)} \cdot
\,C \psi_{\Lambda} \hspace{8mm}
\mbox{for all\/ }C \in {\cal C},\Lambda \in {\cal S}_{\lambda} $$
\hspace*{5mm}where\/ $q_{\qs k}$ is defined by\/
{\em(\ref{q_def})} for $\Lambda$.

\noindent
c) Actions of the elements\/ $X_{k}$ and\/ $\iota(X_{k})$ of\/ $\ASA$ in the
subspace\/ $V_\lambda\otimes_{\,\FF}\KK\subset M_{\chi_{0}}$ coincide.
\et

\noindent
{\it Proof.}
Denote by $V$ the subspace in $\qSprime$ spanned by the elements
$C \psi_{\Lambda}$ where $C \in {\cal C}, \Lambda \in {\cal S}_{\lambda}$.
In the proof of Theorem \ref{fusion_theorem} the element $\psi_{\Lambda}$
is expanded into a sum with leading term $T_{w_{\Lambda}}$.
Thus the elements $C \psi_{\Lambda}$ where $C \in {\cal C}$ and
$\Lambda \in {\cal S}_{\lambda}$ are linearly independent over $\FF$
and form a basis in the space $V$.
We will prove that $V = V_{\lambda}$.

Due to (\ref{eqnA}) every element $\psi_{\Lambda}$ can be
expressed as $\mu_{\Lambda} \,\psi_{\Lambda^{r}}$ for some
invertible $\mu_{\Lambda} \in \qSprime$. So we have
the inclusion $V \subseteq V_{\lambda}$\qs.
The opposite inclusion $V \supseteq V_{\lambda}$ will be established
once we have proved that the action of the
algebra $\qSa$ in $\qSprime$ via left multiplication
preserves the subspace $V$.
The action of the elements $C_{1},\ldots,C_{n}$ preserves the subspace in $V$
\begin{equation}
\label{*}
\underset{C \in {\cal C}}{\oplus} ~\FF\cdot C \psi_{\Lambda}
\end{equation}
for every $\Lambda \in {\cal S}_{\lambda}$,
and thus preserves the space $V$ itself.
We will use Proposition \ref{psi_formulae} to verify that
the action of the elements $T_{1},\ldots,T_{n-1}$ preserves $V$.
Let a standard tableau $\Lambda \in {\cal S}_{\lambda}$ and an index
$k \in \{ 1,\ldots,n-1 \}$ be fixed.
There are three cases to consider.

i) Suppose that \mbox{$s_{k} \cdot \Lambda$} is standard and that
$s_{k}w_{\Lambda} \succ w_{\Lambda}$.
Then Proposition \ref{psi_formulae}(a) gives the equality
$\psi_{k}(q_{\qs k},q_{\qs k+1}) \cdot \psi_{\Lambda} \,=\,
\psi_{s_{k} \cdot \Lambda}$.
Rewriting this equality in the form
\begin{equation}
\label{T1}
T_{k} \cdot \psi_{\Lambda} \;=~ 
\psi_{s_{k} \cdot \Lambda}
\,-\,\left(
\frac{\varepsilon}{q_{\qs k}^{-1}q_{\qs k+1} - 1}
+ \frac{\varepsilon}{q_{\qs k}q_{\qs k+1} - 1} \,C_{k}C_{k+1}
\right)
\psi_{\Lambda}
\end{equation}
explicitly describes the action of the generator $T_{k}$ on the basis
element $\psi_{\Lambda}$.

ii) Next, suppose that \mbox{$s_{k} \cdot \Lambda$} is standard but that
$s_{k}w_{\Lambda} \prec w_{\Lambda}$.
Then Proposition \ref{psi_formulae}(b) yields
$\psi_{k}(q_{\qs k+1},q_{\qs k}) \cdot \psi_{s_{k} \cdot \Lambda} \,=\,
\psi_{\Lambda}$.
Since $\qseq \in {\cal Y}$ then the pair $(q_{\qs k+1},q_{\qs k})$
does not satisfy the condition (\ref{idem_con})\qs; consequently
the element $\psi_{k}(q_{\qs k+1},q_{\qs k})$ is invertible.
Let us multiply the above equality on the left by the element
$\psi_{k}(q_{\qs k},q_{\qs k+1})$. Put
$$
\beta_{k}=
1\,-\,\varepsilon^{2} \!\cdot \left(
\frac{q_{\qs k+1}^{-1}q_{\qs k}}{(q_{\qs k+1}^{-1}q_{\qs k} - 1)^{2}} +
\frac{q_{\qs k}q_{\qs k+1}}{(q_{\qs k}q_{\qs k+1} - 1)^{2}} \right) \in {\FF}.
$$
Then the first relation in Lemma \ref{psi_rel} gives
$
\psi_{k}(q_{\qs k},q_{\qs k+1}) \cdot \psi_{\Lambda}\,=\,
\beta_{k}\,\psi_{s_{k} \cdot \Lambda}\,.
$
Rewriting this equality gives the explicit action of the element $T_{k}$ on the
basis element $\psi_{\Lambda}$ as
\begin{equation}
\label{T2}
T_{k} \cdot \psi_{\Lambda} \;=~
\beta_{k}\,\psi_{s_{k} \cdot \Lambda}
\,-\,\left(
\frac{\varepsilon}{q_{\qs k}^{-1}q_{\qs k+1} - 1}
+ \frac{\varepsilon}{q_{\qs k}q_{\qs k+1} - 1} \,C_{k}C_{k+1}
\right)
\psi_{\Lambda}\,.
\end{equation}

iii) Finally, consider the remaining case where the tableau
\mbox{$s_{k} \cdot \Lambda$} is not standard.
Then Proposition \ref{psi_formulae}(c) gives
$\psi_{k}(q_{\qs k},q_{\qs k+1}) \cdot \psi_{\Lambda} \,=\, 0$.
This equality can be rewritten as
\begin{equation}
\label{T3}
T_{k} \cdot \psi_{\Lambda} \;=\; - \left(
\frac{\varepsilon}{q_{\qs k}^{-1}q_{\qs k+1} - 1}
+ \frac{\varepsilon}{q_{\qs k}q_{\qs k+1} - 1} \,C_{k}C_{k+1}
\right)
\psi_{\Lambda}\,.
\end{equation}

In every instance, the action of the generator $T_{k}$ takes the
element $\psi_{\Lambda}$ to some element within the subspace $V$.
This action can be extended to the full basis in $V$ using the relations
(\ref{sergeev_def3}).
Thus left multiplication in the algebra $\qSprime$ by each generator $T_{k}$
preserves the subspace $V$. 
Then $V_{\lambda} = \qSprime \cdot \psi_{\Lambda^{r}} \subseteq V$
and the part (a) of Theorem \ref{V_theorem} is verified.

Consider the action of the affine generators $X_{1},\ldots,X_{n}$
on the subspace $V_{\lambda}\otimes_{\,\FF}\KK$ in $M_{\chi_{0}}$.
Let us fix the standard tableau $\Lambda \in {\cal S}_{\lambda}$.
Given any $\xseq \in {\cal X}$, we have the equality
$$
\psi_{\Lambda}\xseq \;=\; \pi_{\chi}(\Phi_{w_{\Lambda}})(1)
$$
where $\chi$ is the generic character defined through (\ref{char_def})
by $x_{1},\ldots,x_{n}$.
At the special point $\qseq$ Corollary \ref{scalar_action} then gives
$$
\pi_{\chi_{0}}(X_{k})(\psi_{\Lambda}) ~=~
(w_{\Lambda} \cdot \chi_{0})(X_{k}) \cdot \psi_{\Lambda}
$$
for each $k=1,\ldots,n$.
But using the relations (\ref{affine_rels2}) we obtain
$$
X_{k} \,C \;=\; C \,X_{k}^{-\nu_{C}(k)}, \quad C \in {\cal C}.
$$
Combining these results, for any index $k$ the action of $X_{k}$
on the vector $C \psi_{\Lambda}$ is given by
$$
\pi_{\chi_{0}}(X_{k})(C \psi_{\Lambda}) ~=~
(w_{\Lambda} \cdot \chi_{0})(X_{k}^{-\nu_{C}(k)}) \cdot C \psi_{\Lambda}\;;
\quad C \in {\cal C},\Lambda \in {\cal S}_{\lambda}\,.
$$
The equality stated in Theorem \ref{V_theorem}(b) now follows from the
definition of the character $\chi_{0}$.
In particular, the action of $X_{1},\ldots,X_{n}$ in
$M_{\chi_{0}}$ preserves the subspace $V_{\lambda}\otimes_{\,\FF}\KK$.
Since $\Lambda(1,1) = 1$ for any standard tableau
$\Lambda \in {\cal S}_{\lambda}$, we have $q_{\qs 1}=1$
and the result in (b) gives
$$
\pi_{\chi_{0}}(X_{1})(C \psi_{\Lambda}) ~=~ C \psi_{\Lambda} \;; \quad
C \in {\cal C},\ \Lambda \in {\cal S}_{\lambda}\,.
$$
Thus the actions of $X_{1}$ and $\iota(X_{1}) = 1$ in
$V_{\lambda}\otimes_{\,\FF}\KK$ coincide\qs. Theorem \ref{V_theorem}(c)
now follows from Proposition \ref{Murphy_hom}.\hspace{3mm}\BOX

\smallskip\noindent
The basis in $V_{\lambda}$ described in Theorem \ref{V_theorem}(a)
is an analogue of the Young basis \cite{DJ} in an
irreducible $H_{n}(q)$-module.
Theorem \ref{V_theorem}(a) also shows that
\mbox{$\dim V_{\lambda} \,=\, 2^{n} \cdot m_{\lambda}$}
where $m_{\lambda}$ is the number of standard shifted tableaux with
shape $\lambda$.
The integer $m_{\lambda}$ can be computed by an analogue of the
hook formula \cite[Example I.5.2]{Mac} for the number $n_{\omega}$
of standard Young tableaux with shape $\omega$.
This analogue was first
found by Morris \mbox{\cite[Theorem 2.1]{Mor}\qs;} see also
\cite[\mbox{Example III.7.8}]{Mac}.

For the remainder of this paper, we will consider $V_{\lambda}$ as a
$\qSprime$-module. The vector space $V_\lambda$ over $\FF$
has a natural decomposition into the direct sum of the subspaces (\ref{*}).
For each $i=1,\ldots,\ell_{\lambda}$
define a linear operator $\rho_i$ in $V_\lambda$ by
$$
\rho_i\,:\ C\,\psi_\Lambda\,\mapsto\,C\,C_l\,\psi_\Lambda\,,
\quad
l=\Lambda(i,i)
$$
where $\Lambda\in{\cal S}_{\lambda}$. 
Then we have the following easy observation; cf.\  
\text{\cite[Theorem 8.3]{Naz}.}

\bp
\label{commutant}
The operators $\rho_1,\ldots,\rho_{\ell_\lambda}$ commute
with the action of $\qSprime$ in $V_\lambda$\qs.
\ep

\noindent
{\it Proof.}
By definition the operators $\rho_1,\ldots,\rho_{\ell_\lambda}$
commute with the action of $C_1,\ldots,C_n\in\qSprime$ in $V_\lambda$\qs.
Now fix an index $i\in\{1,\ldots,\ell_\lambda\}$ and a tableau
$\Lambda\in\cal S_\lambda$. Put $l=\Lambda(i,i)$.
It suffices to prove that in $V_\lambda$
$$
T_k\cdot \rho_i(\psi_\Lambda)=
T_k\,C_l\,\psi_\Lambda=
\rho_i(T_k\,\psi_\Lambda)
$$
for any $k\in\{1,\ldots,n-1\}$. For $k\neq l,l-1$ the required second
equality immediately follows from the third relation in (\ref{sergeev_def3})
and from (\ref{T1}),(\ref{T2}),(\ref{T3}).

Suppose that the tableau $s_{k}\cdot\Lambda$
is standard and that
$s_{k} w_{\Lambda} \succ w_{\Lambda}$. The cases when 
$s_{k}\cdot\Lambda$ is not standard or it is but 
$s_{k}w_{\Lambda} \prec w_{\Lambda}$ can be treated similarly,
see the beginning of the proof of Theorem \ref{V_theorem}.
Let $l=k$, then $q_{\hskip1pt k}=1$ by (\ref{q_def}). Write
$q_{\hskip1pt k+1}=\delta$\qs.
Then by the first relation in (\ref{sergeev_def3}) and by (\ref{T1}) 
we obtain the equalities
\begin{eqnarray*}
T_k\,C_k\,\psi_\Lambda=C_{k+1}\,T_k\,\psi_\Lambda
\hskip-8pt&=&\hskip-7pt
C_{k+1}\,
\psi_{s_k\cdot\Lambda}
-
C_{k+1}
\left(
\frac{\varepsilon}{\delta-1}
+
\frac{\varepsilon}{\delta-1} \,C_{k}C_{k+1}
\right)
\psi_{\Lambda}
\\
\hskip-8pt&=&\hskip-7pt
C_{k+1}\,
\psi_{s_k\cdot\Lambda}
-
\left(
\frac{\varepsilon}{\delta-1}
+
\frac{\varepsilon}{\delta-1} \,C_{k}C_{k+1}
\right)C_k\,
\psi_{\Lambda}
=\rho_i(T_k\,\psi_\Lambda)\,.
\end{eqnarray*}

Now let $l=k+1$, then $q_{\hskip1pt k+1}=1$. Write
$q_{\hskip1pt k}=\delta$.
Then by the second relation in (\ref{sergeev_def3}) and by (\ref{T1})
we obtain the equalities 
\begin{eqnarray*}
T_k\,C_{k+1}\,\psi_\Lambda
\hskip-8pt&=&\hskip-7pt
C_k\,
\psi_{s_k\cdot\Lambda}
-
C_k
\left(
\frac{\varepsilon}{\delta^{\,-1}-1}
+
\frac{\varepsilon}{\delta-1} \,C_{k}C_{k+1}
\right)
\psi_{\Lambda}
+(\varepsilon\,C_{k+1}-\varepsilon\,C_k)\,\psi_\Lambda
\\
\hskip-8pt&=&\hskip-7pt
C_k\,
\psi_{s_k\cdot\Lambda}
\,-\,
\left(
\frac{\varepsilon}{\delta^{\,-1}-1}
+
\frac{\varepsilon}{\delta-1} \,C_{k}C_{k+1}
\right)C_{k+1}\,
\psi_{\Lambda}
=\rho_i(T_k\,\psi_\Lambda)\,.
\hspace{3mm}\BOX
\end{eqnarray*}

\smallskip\noindent
Put $d_{\lambda} = 0$ if the number $\ell_\lambda$ is even,
and put $d_{\lambda} = 1$ if $\ell_\lambda$ is odd.
The assignment $C_i\mapsto \rho_i$ defines an action in $V_\lambda$
of the Clifford algebra $Z_\lambda$ over ${\mathbb C}$
with $\ell_\lambda$ generators. This algebra
has a natural \Z -gradation: each of its generators 
$C_1,\ldots,C_{\ell_\lambda}$ is odd.
Take any minimal even idempotent in this algebra.
For example, take any of the $2^{[\,\ell_{\lambda}/2\,]}$
pairwise orthogonal idempotents
$$
2^{-[\,\ell_{\lambda}/2\,]}\cdot
\left(\,
1\pm C_1C_2\,\sqrt{-1}
\,\right)
\dots
\left(\,
1\pm C_{\ell_\lambda-d_\lambda-1}\,C_{\ell_\lambda-d_\lambda}\sqrt{-1}
\,\right).\hskip-10pt
$$
The left ideal in $Z_\lambda$ generated by this idempotent is
irreducible
under the left multiplication. It is absolutely irreducible if and only if the
number $\ell_\lambda$ is even.
Let $\gamma_{\Lambda^r}$ be the image of this idempotent under
the embedding of $Z_\lambda$ into the Clifford algebra with $n$
generators over ${\mathbb C}$ determined by
$$
C_i\mapsto C_l\,,
\quad
l=\Lambda^r(i,i)\,.
$$ 
Consider the left ideal $U_\lambda$ in the algebra $\qSprime$
generated by the element
$\gamma_{\Lambda^r}\,\psi_{\Lambda^r}$
as a $\qSprime$-module under the left multiplication.
Then we get the following corollary to \text{Proposition \ref{commutant}.} 

\bc
\label{splitting}
The $\qSprime$-module $V_\lambda$ 
is a direct sum of\/ $2^{[\,\ell_{\lambda}/2\,]}$
copies of the module\/ $U_{\lambda}$.
\ec

\noindent
By the definition (\ref{q_def}) we get $q_{\qs l}=1$ if $l=\Lambda(i,i)$.
Denote by $\gamma_{\Lambda^r}^{\,\prime}$ the image of the idempotent
$\gamma_{\Lambda^r}$ under the automorphism of
the Clifford algebra with $n$ generators over ${\mathbb C}$ determined by
$$
C_l\,\mapsto\,C_{w_{\Lambda^r}^{-1}(l)}\,\,;
\ \ \quad
l=1,\ldots,n\,.
$$ 
Then the next proposition follows immediately from Proposition \ref{c_inter}.

\bp
\label{gamma_inter}
We have the equality\/
$
\gamma_{\Lambda^r}\hskip1pt\psi_{\Lambda^r}=
\psi_{\Lambda^r}\hskip1pt\gamma_{\Lambda^r}^{\,\prime}
$
in the algebra $\qSprime$.
\ep

\noindent
We now describe the {\it centre\/}
$\operatorname{Z}(\qSa)$ of the \Z-graded algebra $\qSa$.
This is the collection of elements $Z\in \qSa$ such that
the supercommutator $[\qs Z\hskip1pt,X\qs ]$ vanishes for all $X \in \qSa$.
Note that the even component of $\operatorname{Z}(\qSa)$
coincides with the even component of the centre in the usual sense.
We use the Jucys-Murphy elements $J_1,\ldots,J_n\in\qSa$
defined by (\ref{Jucys_def}).

\bp
\label{Sergeev_centre}
a) The centre\/ $\operatorname{Z}(\qSa)$
of the\/ \Z -graded algebra\/ $\qSa$ consists precisely of
the symmetric polynomials in the elements\/
$J_{1} + J_{1}^{-1},\ldots,J_{n} + J_{n}^{-1}$\qs.

\noindent
b) The dimension of\/ $\operatorname{Z}(\qSa)$
over\/ ${\mathbb C}(q)$ coincides with the
number of strict partitions\/ $\lambda \stp n$.
\ep

\noindent
{\it Proof.}
Propositions \ref{max_comm}(b) and \ref{Murphy_hom}
imply that any symmetric polynomial in
the elements $J_{1} + J_{1}^{-1},\ldots,J_{n} + J_{n}^{-1}$
is central (in the usual sense) for the algebra $\qSa$.
Furthermore, these central elements are even since the Jucys-Murphy elements
$J_{1},\ldots,J_{n}$ are homogeneous with degree zero in the \Z -grading
of $\qSa$.
We have to prove that all these symmetric polynomials exhaust
the centre $\operatorname{Z}(\qSa)$, cf. \cite[Section 3]{Jucy2}.

Let $d$ denote the number of strict partitions of $n$.
The algebra $\Sa$ can be obtained from the generic
algebra $\qSa$ via the specialization $q \rightarrow 1$.
The dimension of the centre does not decrease under this specialization\qs.
But the centre of the \Z -graded algebra $\Sa$ has dimension~$d$
\cite[Lemma 6]{Serg}.
Thus the integer $d$ provides an upper bound on the dimension
of $\operatorname{Z}(\qSa)$.
It remains to show that the collection of symmetric polynomials in
$J_{1} + J_{1}^{-1},\ldots,J_{n} + J_{n}^{-1}$
contains at least $d$ linearly independent elements over $\CCq$.

By the part (a) of Theorem \ref{V_theorem} the elements $C \psi_{\Lambda}$
with $\Lambda \in {\cal S}_{\lambda}$ and $C \in {\cal C}$
form a $\FF$-basis in $V_\lambda$. By the part (b) of the same theorem we have
$$
(J_{k} + J_{k}^{-1}) \cdot C \psi_{\Lambda} \;=\;
(q_{\qs k} + q_{\qs k}^{-1}) \cdot C \psi_{\Lambda} \;=\;
\frac{2}{q+q^{-1}} \left( \,q^{2(j-i)+1} + q^{-2(j-i)-1} \right)
\cdot C \psi_{\Lambda}
$$
where $k = \Lambda(i,j)$\qs; see (\ref{scal_sub}) and (\ref{q_def}).
So any symmetric polynomial in
$J_{1} + J_{1}^{-1},\ldots,J_{n} + J_{n}^{-1}$
acts in $V_\lambda$ by a certain scalar from $\CCq$.
The collection of these scalars as we range over all symmetric polynomials,
determines the partition $\lambda$ uniquely.
Indeed, the generating function
$(\,t+J_{1}+J_{1}^{-1})\,\cdots\,(\,t+J_{n}+J_{n}^{-1})$
in $t$ for the elementary symmetric polynomials has in $V_\lambda$
the eigenvalue
$$
2^{n} \cdot \prod_{k=1}^{n}
\left( \,\frac{t}{2} + \frac{q^{2(j-i)+1}+q^{-2(j-i)-1}}{q+q^{-1}} \,\right).
$$
It allows us to recover the unordered collection
of the contents $j-i$ for the shifted Young diagram of $\lambda$\qs;
this collection determines the partition $\lambda$ uniquely.
Therefore the elementary symmetric polynomials in
$J_{1} + J_{1}^{-1},\ldots,J_{n} + J_{n}^{-1}$
generate at least $d$ linearly independent elements of
$\operatorname{Z}(\qSa)$.\hspace{3mm}\BOX

\bt
\label{V_results}
The module\/ $U_\lambda$ over the \Z-graded algebra\/ $\qSprime$

a) is irreducible\/,

b) remains irreducible on passing to any extension of the field\/ $\FF$,

c) is absolutely irreducible if and only if the number\/ 
$\ell_\lambda$ is even\/.
\et

\noindent
{\it Proof.}
Let $U_\lambda^{\hspace{0.45mm}\text{\bf-}}$
be a non-zero irreducible submodule in the $\qSdash$-module
$U_\lambda\otimes_{\,\FF}\KK$.
In the proof of Proposition \ref{Sergeev_centre}
we demonstrated that the central elements of the \mbox{\Z-graded}
algebra $\qSa$ act in $V_{\lambda}$ by certain scalars.
These central elements act in  $U_\lambda^{\hspace{0.45mm}\text{\bf-}}$
by the same scalars. These scalars distinguish the modules
 $U_\lambda^{\hspace{0.45mm}\text{\bf-}}$
for different strict partitions $\lambda\stpscript n$.

The $\qSprime$-module $U_\lambda$ with $d_\lambda=1$
cannot be absolutely irreducible by Proposition~\ref{commutant}.
But Theorem \ref{V_theorem}(a) and Corollary \ref{splitting}
along with the classical result of \cite{Sch} show that 
$$
\sum_{\lambda \stpscript n}\,2^{\,-d_{\lambda}} \cdot
\bigl(\,\operatorname{dim}U_{\lambda}\,\bigr)^{2} \,=\,
\sum_{\lambda \stpscript n}
\bigl(\,m_{\lambda} \cdot 2^{\,n-(\ell_{\lambda}/2)}\,\bigr)^{2} \,=\;
2^{n}\cdot n \hspace{1pt}!
$$
which is exactly the dimension of the algebra $\qSprime$ over the field $\FF$,
see Proposition \ref{dimension}.
So if  $U_\lambda^{\hspace{0.45mm}\text{\bf-}}$ is a proper submodule in
$U_\lambda\otimes_{\,\FF}\KK$
or if the $\qSprime$-module $U_\lambda$ with $d_\lambda=0$ is not
absolutely irreducible, we get
a contradiction with Propositions \ref{semisimple} and \ref{Sergeev_centre}(b).
\hspace{3mm}\BOX

\smallskip\noindent
In the course of the proof of Theorem \ref{V_results}
we also established the following fact.

\bc
\label{complete_set}
The modules\/ $U_{\lambda}$ ranging over all strict partitions\/
$\lambda\stp n$
form a complete set of irreducible pairwise non-equivalent\/
$\qSprime$-modules.
\ec

\smallskip\noindent
Thus we have demonstrated that $\FF$ is a splitting field for
the semisimple $\CCq$-algebra $\qSa$.
\newcommand{\BOOK}[6]{\bibitem[{\bf #6}]{#1}{\sc #2}, {\it #3} (#4)#5.}
\newcommand{\JPAPER}[8]{\bibitem[{\bf #8}]{#1}{\sc #2}, `#3',
{\it #4} #5 (#6) #7.}
\newcommand{\JPAPERS}[9]{\bibitem[{\bf #9}]{#1}{\sc #2}, `#3', {\it #4} #5 #6,
#7 #8.}
%

%
%
\begin{it}
Mathematics Department
\hfill
E-mail\::\hspace{1mm}
{\em arj4@york.ac.uk}

\noindent
University of York
\hfill
{\em mln1@york.ac.uk}

\noindent
York YO1 5DD

\noindent
England

\end{it}

\end{document}